\documentclass[12pt]{amsart}
\usepackage{amscd}
\usepackage{amssymb}

\allowdisplaybreaks[1]

\newtheorem{thm}{Th\'eor\`eme}[section]
\newtheorem{rema}[thm]{Remarque}

\newtheorem{lem}[thm]{Lemme}
\newtheorem{prop}[thm]{Proposition}

\newtheorem{defn}[thm]{D\'efinition}
\theoremstyle{definition}

\numberwithin{equation}{section}

\newcommand{\R}{\mathbb R}

\begin{document}

\date{30/06/05}

\title[Cohomologie vectorielle]{Cohomologie de Chevalley des
graphes vectoriels}
\author[W. Aloulou, D. Arnal et R. Chatbouri]{Walid Aloulou, Didier Arnal et 
Ridha Chatbouri}
\address{
D\'epartement de Math\'ematiques\\
Unit\'e de Recherche Physique Math\'ematique\\
Fa\-cult\'e des Sciences de Monastir\\
Avenue de l'environnement\\
5019 Monastir\\
Tuni\-sie} \email{Walid.Aloulou@ipeim.rnu.tn}

\address{
Institut de Math\'ematiques de Bourgogne\\
UMR CNRS 5584\\
Universit\'e de Bourgogne\\
U.F.R. Sciences et Techniques
B.P. 47870\\
F-21078 Dijon Cedex\\France} \email{Didier.Arnal@u-bourgogne.fr}

\address{
D\'epartement de Math\'ematiques\\
Unit\'e de Recherche Physique Math\'ematique\\
Fa\-cult\'e des Sciences de Monastir\\
Avenue de l'environnement\\
5019 Monastir\\
Tuni\-sie} \email{Ridha.Chatbouri@ipeim.rnu.tn}

\thanks{
Ce travail a \'et\'e effectu\'e dans le cadre de l'accord CMCU 03 S 1502. 
W. Aloulou et R. Chatbouri remercient l'universit\'e de Bourgogne pour l'accueil
dont ils ont b\'en\'efici\'e au cours de leurs s\'ejours, D. Arnal remercie la 
facult\'e des Sciences de Monastir pour l'accueil dont il a b\'en\'efici\'e au 
cours de ses s\'ejours.}


\begin{abstract}
L'espace des champs de vecteurs augment\'e des fonctions $C^\infty$ sur 
${\mathbb R}^d$ est une sous alg\`ebre de Lie de l'alg\`ebre de Lie (gradu\'ee) 
de l'espace $T_{poly}({\mathbb R}^d)$ des champs de tenseurs sur ${\mathbb R}^d$ 
muni du crochet de Schouten.

Dans cet article, on calcule la cohomologie des repr\'esentations adjointes de 
cette sous alg\`ebre de Lie, en se restreignant \`a des cocha\^\i nes d\'efinies 
par des graphes de Kontsevich a\'eriens comme dans \cite {[AGM]}. On retrouve 
les r\'esultats bien connus de \cite{[GF]} et \cite{[DWL]}.

\end{abstract}


\maketitle


\

\section{Introduction}\label{sec1}

\

Notons $T_{poly}({\mathbb R}^d)$ l'espace des tenseurs compl\`etement
antisym\'etriques sur ${\mathbb R}^d$. Cet espace, muni du crochet de
Schouten et de la graduation $deg(\alpha)=k-1$ si $\alpha$ est un
$k$-tenseur est une alg\`ebre de Lie gradu\'ee.

Cette alg\`ebre de Lie contient une sous-alg\`ebre de Lie int\'eressante : 
$Vect({\mathbb R}^d)$, espace des tenseurs de degr\'e n\'egatif ou 
nul, c'est \`a dire l'alg\`ebre de Lie des champs de vecteurs $\xi$ augment\'ee 
de l'espace des fonctions $C^\infty$ $f$ muni du crochet usuel des champs de 
vecteurs \'etendu aux fonctions par:
$$
[\xi,f]=-[f,\xi]=\xi f,~~\hbox{et}\quad [f_1,f_2]=0.
$$

La repr\'esentation adjointe fait de $T_{poly}({\mathbb R}^d)$ un
$Vect({\mathbb R}^d)$-module. Dans cet article, on va calculer des groupes de
cohomologies de ce module.

Plus pr\'ecis\'ement, nous consid\'erons ici des cocha\^\i nes $\varphi$ 
d\'efinies comme dans \cite{[AGM]} \`a partir d'une combinaison lin\'eaire de 
graphes a\'eriens de Kontsevich \cite{[K]}. Dans le cas de 
$Vect({\mathbb R}^d)$, au plus une ar\^ete part d'un sommet de chacun de ces 
graphes. L'op\'erateur de cohomologie peut alors \^etre d\'efini sur l'espace 
des graphes, il correspond \`a une suite d'`\'eclatement' des sommets.

Les cohomologies de Chevalley des champs de vecteurs ont \'et\'e calcul\'ees par
plu\-sieurs auteurs. En particulier dans \cite{[DWL]}, la cohomologie \`a 
valeurs dans les formes est d\'etermin\'ee explicitement. Dans cet article, 
notre restriction aux graphes purement a\'eriens nous permet d'adapter la preuve 
de \cite{[DWL]}. On reformule en particulier la d\'efinition de l'homotopie et 
on en d\'eduit la cohomologie des graphes a\'eriens. On retrouve le m\^eme 
r\'esultat, \`a savoir que cette cohomologie est donn\'ee par les roues 
impaires.

Rappelons que dans le probl\`eme de la construction d'une formalit\'e sur 
${\mathbb R}^d$ au moyen de graphes, les cohomologies qui apparaissent sont 
celles de Hochschild, de Chevalley et de Chevalley-Harrison. La premi\`ere a
\'et\'e calcul\'ee dans \cite{[AM]}, la troisi\`eme est nulle d'apr\`es 
\cite{[GH]}.


\

\section{Notations et d\'efinitions}

\

Dans cet article, on consid\`ere l'espace $T_{poly}(\R^d)$ des
tenseurs compl\`etement antisym\'etriques sur $\R^d$. Si
${\mathcal X}$ d\'esigne l'espace des champs de vecteurs $\xi$ sur
$\R^d$, on construit $Tens(\R^d)$ comme l'alg\`ebre associative
libre sur $C^\infty(\R^d)$ engendr\'ee par les champs constants
$\partial_1,\dots,\partial_d$. $T_{poly}(\R^d)$ est le quotient de
$Tens(\R^d)$ par l'id\'eal engendr\'e par
$\{\xi\otimes\eta-\eta\otimes\xi, ~~\xi\in{\mathcal
X},~~\eta\in{\mathcal X}\}$. $T_{poly}(\R^d)$ est muni d'un
produit associatif $\wedge$, tout \'el\'ement de $T_{poly}(\R^d)$
est une somme de produits de la forme
$\xi_1\wedge\dots\wedge\xi_k$ et de fonctions $f$. On peut aussi
\'ecrire tout tenseur $\alpha$ de $T_{poly}(\R^d)$ de fa\c con
unique sous la forme
$$
\alpha=\sum_{k=0}^K\sum_{i_1,\dots,i_k}\alpha_{(k)}^{i_1\dots i_k}
\partial_{i_1}\wedge\dots\wedge\partial_{i_k}.
$$
(On prend la convention que les coordonn\'ees
$\alpha_{(k)}^{i_1\dots i_k}$ sont des fonctions $C^\infty$ et
sont compl\`etement antisym\'etriques en $i_1,\dots, i_k$. On a
donc aussi
$$
\alpha=\sum_{k=0}^K\sum_{i_1<\dots<i_k}~k!~\alpha_{(k)}^{i_1\dots i_k}
\partial_{i_1}\wedge\dots\wedge\partial_{i_k}.
$$

On place sur $\R^d$ la connexion plate triviale $\nabla$, c'est \`a dire la 
connexion pour la structure riemannienne usuelle de $\R^d$. On a donc
$$
\nabla_\xi f=\xi f,~~\nabla_\xi\eta=\nabla_\xi(\sum_i\eta^i\partial_i)=\sum_i
(\xi \eta^i)\partial_i\quad(f\in C^\infty(\R^d),~\xi,~\eta\in{\mathcal X}).
$$
Il y a un prolongement unique de $\nabla_\xi$ en une d\'erivation de
$T_{poly}(\R^d)$. On impose
$$
\nabla_\xi(\alpha\wedge\beta)=(\nabla_\xi\alpha)\wedge\beta+\alpha\wedge
(\nabla_\xi\beta),
$$
on obtient une solution et une seule d\'efinie par
$$
\nabla_\xi(\eta_1\wedge\dots\wedge\eta_\ell)=\sum_{j=1}^\ell(-1)^{j-1}(
\nabla_\xi\eta_j)\wedge\eta_1\wedge\dots\widehat{\eta_j}\dots\wedge\eta_\ell
$$
ou par
$$
\nabla_\xi\alpha=\sum_{i_1,\dots,i_k}\left(\xi\alpha^{i_1\dots i_k}\right)
\partial_{i_1}\wedge\dots\wedge\partial_{i_k}.
$$

Maintenant, on consid\`ere $T_{poly}(\R^d)$ comme une alg\`ebre gradu\'ee par 
$|\alpha|=k$ si $\alpha$ est un $k$-tenseur et on utilise syst\'ematiquement la 
r\`egle de Koszul. Si $\xi\in{\mathcal X}$, $\nabla_\xi$ est une d\'erivation de 
degr\'e 0, le produit $\wedge$ est aussi de degr\'e 0, les formules
ci-dessus sont donc coh\'erentes avec cette r\`egle. Cependant l'application 
$\nabla:\xi\mapsto\nabla_\xi$ est maintenant homog\`ene de degr\'e -1 de 
${\mathcal X}$ vers l'espace $Der(T_{poly}(\R^d))$ des d\'erivations de
$T_{poly}(\R^d)$. On veut la prolonger comme une d\'erivation. 

On cherche donc $\nabla:T_{poly}(\R^d)\longrightarrow Der(T_{poly}(\R^d))$ qui 
la prolonge et telle que
$$
\nabla_{\alpha\wedge\beta}(\gamma)=(-1)^{|\alpha|}\alpha\wedge(\nabla_\beta
\gamma)+(-1)^{|\beta||\gamma|}(\nabla_\alpha\gamma)\wedge\beta.
$$

\

\begin{lem} {\rm (Le prolongement)}

\

Un tel prolongement existe et il est unique, il est d\'efini par
$\nabla_f=0$ si $f\in C^\infty(\R^d)$ et soit
$$
\nabla_{\xi_1\wedge\dots\wedge\xi_k}(\eta_1\wedge\dots\wedge\eta_\ell)=
\sum_{i=1}^k\sum_{j=1}^\ell~(-1)^{i+j}\xi_1\wedge\dots\widehat{\xi_i}\dots\wedge
\xi_k\wedge(\nabla_{\xi_i}\eta_j)\wedge\eta_1\wedge\dots\widehat{\eta_j}\dots
\wedge\eta_\ell
$$
soit
$$
\nabla_\alpha\beta=\sum_{r=1}^k~(-1)^{r-1}~\sum_{i_1,\dots,i_k \atop j_1,\dots
,j_\ell}\alpha^{i_1\dots i_k}\left(\partial_{i_r}\beta^{j_1\dots j_\ell}\right)
\partial_{i_1}\wedge\dots\widehat{\partial_{i_r}}\dots\wedge\partial_{i_k}\wedge
\partial_{j_1}\wedge\dots\wedge\partial_{j_\ell}.
$$
\end{lem}

\

Remarquons que l'on retrouve l'op\'eration not\'ee $\bullet$ dans \cite{[AMM]}
$$
\nabla_\alpha\beta=\alpha\bullet\beta.
$$\\

\vskip 0.3cm
\noindent
{\bf Preuve}

Soit $\xi$ un champ de vecteur et $f\in C^\infty(\R^d)$. On doit
avoir $\nabla_{(f\xi)}=\nabla_{f\wedge\xi}$ ou, pour tout
$k$-tenseur $\alpha$,
$$
f\nabla_\xi\alpha=\nabla_{(f\xi)}(\alpha)=f(\nabla_\xi\alpha)+(-1)^k
(\nabla_f\alpha)\wedge\xi.
$$
$\nabla_f\alpha$ est un $k-1$-tenseur. S'il n'est pas nul, on peut
choisir $\xi$ tel que $(\nabla_f\alpha)\wedge\xi\neq0$, ce qui est
absurde.

On montre ensuite par r\'ecurrence sur $k$ que
$\nabla_{\xi_1\wedge\dots\wedge
\xi_k}(\eta_1\wedge\dots\wedge\eta_\ell)$ ne peut \^etre que ce
qui est annonc\'e. Enfin que l'application $\nabla$ ainsi
d\'efinie a bien les propri\'et\'es demand\'ees. La formule
donnant $\nabla_\alpha\beta$ est une
cons\'equence imm\'ediate de la premi\`ere formule.\\

Maintenant on a par construction
$$\nabla_\xi\eta-\nabla_\eta\xi=[\xi,\eta],\quad(\xi,~\eta\in{\mathcal X}).$$
Traditionnellement, pour \'etendre le crochet des champs de
vecteurs en le crochet de Schouten, d\'efini sur $T_{poly}(\R^d)$,
on choisit la graduation $deg(\alpha)=k-1$ si $\alpha$ est un
$k$-tenseur, le crochet des champs de vecteurs est
antisym\'etrique, de degr\'e 0, il se prolonge d'une fa\c con
unique en un crochet antisym\'etrique toujours de degr\'e 0 sur
$T_{poly}(\R^d)$ qui est une d\'erivation '\`a droite' c'est \`a
dire en un crochet tel que
\begin{align*}
[\alpha,\beta\wedge\gamma]&=[\alpha,\beta]\wedge\gamma+(-1)^{(deg(\beta)-1)
deg(\alpha)}\beta\wedge[\alpha,\gamma]\cr
[\beta,\alpha]&=-(-1)^{deg(\alpha)deg(\beta)}[\alpha,\beta]
\end{align*}
(en effet, $deg(\wedge)=-1$ maintenant). Ce prolongement unique
est donn\'e par $[f,g]=0$ si $f$ et $g$ sont des fonctions,
$[\xi,f]=-[f,\xi]=\xi f$ si $\xi$ est un champ de vecteurs et $f$
une fonction et par
$$[\xi_1\wedge\dots\wedge\xi_k,\eta_1\wedge\dots\wedge\eta_\ell]=\sum_{i=1}^k
\sum_{j=1}^\ell(-1)^{k-i+j-1}\xi_1\wedge\dots\widehat{\xi_i}\dots\wedge\xi_k
\wedge[\xi_i,\eta_j]\wedge\eta_1\wedge\dots\widehat{\eta_j}\dots\wedge\eta_\ell
.$$ Il v\'erifie la relation
$$[\alpha\wedge\beta,\gamma]=\alpha\wedge[\beta,\gamma]+(-1)^{(deg(\beta)+1)
deg(\gamma)}[\alpha,\gamma]\wedge\beta.$$

Avec nos notations, on peut aussi d\'efinir le crochet de Schouten
par:
$$[\alpha,\beta]=(-1)^{deg(\alpha)}\nabla_\alpha\beta-(-1)^{(deg(\alpha)+1)
deg(\beta)}\nabla_\beta\alpha.$$

Cependant, nous gardons ici la graduation $|\alpha|=k$ si $\alpha$
est un $k$-tenseur. Le crochet des champs de vecteurs devient un
produit commutatif. On le prolonge donc comme dans \cite{[K]} en
une op\'eration $Q$ sym\'etrique, de degr\'e -1 v\'erifiant
\begin{align*}
Q(\alpha,\beta\wedge\gamma)&=Q(\alpha,\beta)\wedge\gamma+(-1)^{|\beta|(|
\alpha|-1)}\beta\wedge Q(\alpha,\gamma)\cr
Q(\beta,\alpha)&=(-1)^{|\alpha||\beta|}Q(\alpha,\beta).
\end{align*}

\begin{lem} {\rm (L'op\'erateur $Q$)}

\

Il y a un prolongement et un seul $Q$ du crochet des champs de vecteurs 
v\'erifiant ces relations. Ce prolongement est donn\'e par
$$
Q(\alpha,\beta)=\nabla_\alpha\beta+(-1)^{|\alpha||\beta|}\nabla_\beta\alpha.
$$
Il v\'erifie
$$
Q(\alpha\wedge\beta,\gamma)=(-1)^{|\alpha|}\alpha\wedge Q(\beta,\gamma)+
(-1)^{|\beta||\gamma|}Q(\alpha,\gamma)\wedge\beta.
$$
On a
\begin{align*}
Q(\alpha,\beta)&=\sum_{\begin{smallmatrix} i_1\dots i_k\cr j_1\dots j_\ell
\end{smallmatrix}}~~\Big[\sum_{r=1}^k (-1)^{r-1}\alpha^{i_1\dots i_k}\left(
\partial_{i_r}\beta^{j_1\dots j_\ell}\right)\partial_{i_1}\wedge\dots
\widehat{\partial_{i_r}}\dots\wedge\partial_{i_k}\wedge\partial_{j_1}\wedge
\dots\wedge\partial_{j_\ell}\cr
&\hskip 1cm+\sum_{s=1}^\ell(-1)^{k+s}\beta^{j_1\dots j_\ell}\left(\partial_{j_s}
\alpha^{i_1\dots i_k}\right)\partial_{i_1}\wedge\dots\wedge\partial_{i_k}\wedge
\partial_{j_1}\wedge\dots\widehat{\partial_{j_s}}\dots\wedge\partial_{j_\ell}
\Big].
\end{align*}
et
$$[\alpha,\beta]=(-1)^{deg(\alpha)}Q(\alpha,\beta).$$
\end{lem}

\

\noindent {\bf Preuve}

Supposons que $Q$ soit une telle extension. $Q$ \'etant de degr\'e
-1, $Q(f,g)=0$ si $f$ et $g$ sont des fonctions. Maintenant
$$Q(\xi,f\eta)=Q(\xi,f)\eta+fQ(\xi,\eta)=[\xi,f\eta]=(\xi f)\eta+f[\xi,\eta],$$
donc $Q(\xi,f)=Q(f,\xi)=\xi f=\nabla_\xi f+\nabla_f\xi$.
L'application $\beta \mapsto Q(\xi,\beta)$ \'etant une
d\'erivation, on montre par r\'ecurrence sur $\ell$ que
$$Q(\xi,\eta_1\wedge\dots\wedge\eta_\ell)=\sum_{j=1}^\ell(-1)^{j-1}[\xi,\eta_j]
\wedge\eta_1\wedge\dots\widehat{\eta_j}\dots\wedge\eta_\ell.$$ Par
r\'ecurrence sur $k$, on montre ensuite que
$$Q(\xi_1\wedge\dots\wedge\xi_k,\eta_1\wedge\dots\eta_\ell)=\sum_{
\begin{smallmatrix}i=1\cr j=1\end{smallmatrix}}^{\begin{smallmatrix}i=k\cr
j=\ell\end{smallmatrix}}(-1)^{i+j}\xi_1\wedge\dots\widehat{\xi_i}\dots\wedge
\xi_k\wedge[\xi_i,\eta_j]\wedge\eta_1\wedge\dots\widehat{\eta_j}\dots\wedge
\eta_\ell.$$

Ceci nous dit que si $Q$ existe, elle est unique et que c'est
$$Q(\alpha,\beta)=\nabla_\alpha\beta+(-1)^{|\alpha||\beta|}\nabla_\beta\alpha.$$
Maintenant les propri\'et\'es de $\nabla$ montrent que cette
formule d\'efinit bien une bi\-d\'e\-ri\-va\-tion sym\'etrique de degr\'e
-1, extension du crochet des
champs de vecteurs. Les derni\`eres formules du lemme sont imm\'ediates.\\

Reprenons maintenant le choix de signes donn\'e dans \cite{[AMM]}.
Soit $\sigma$ une permutation de $\{1,\dots,n\}$. On note
$\varepsilon(\sigma)$ sa signature. Si $v_1,\dots,v_n$ sont $n$
vecteurs homog\`enes d'un espace vectoriel $V$ gradu\'e par $deg$,
on note $\varepsilon_{deg(v)}(\sigma)$ la signature de la
permutation que $\sigma$ induit sur les $v$ de degr\'es impairs.
Par construction ces signatures sont des homomorphismes de groupe.
Si $C$ est une application $n$-lin\'eaire sur $V^n$, \`a valeurs
dans un espace vectoriel, on dira que $C$ est sym\'etrique (resp.
antisym\'etrique) si pour toute permutation $\sigma\in S_n$ et
tout $v_1,\dots,v_n$ homog\`enes,
\begin{align*}
C(v_{\sigma(1)},\dots,v_{\sigma(n)})&=\varepsilon_{deg(v)}(\sigma)C(v_1,\dots,
v_n)\cr
\Big(\hbox{respectivement}~~C(v_{\sigma(1)},\dots,v_{\sigma(n)})&=\varepsilon
(\sigma)\varepsilon_{deg(v)}(\sigma)C(v_1,\dots,v_n)\Big).
\end{align*}
Ceci est \'equivalent \`a 
\begin{align*}
C(v_1,\dots,v_{i+1},v_i,\dots,v_n)&=(-1)^{deg(v_i)deg(v_{i+1})}C(v_1,\dots,v_n)
\cr
\big(\hbox{resp.}~~C(v_1,\dots,v_{i+1},v_i,\dots,v_n)&=-(-1)^{deg(v_i)
deg(v_{i+1})}C(v_1,\dots,v_n)\big)
\end{align*}
pour tout $i$.\\

Changeons de graduation sur $V$ et posons $|v|=deg(v)+1$. Posons
$$\eta_v(\sigma)=\varepsilon_{deg(v)}(\sigma)\varepsilon_{|v|}(\sigma)
\varepsilon(\sigma).$$
Si nous nous donnons une application $\tau$ de $V^n$ dans $\{\pm1\}$ telle que 
pour tout $\sigma$ et tout $v_1,\dots,v_n$,
\begin{align}\label{(*)}
\tau(v_{\sigma(1)},\dots,v_{\sigma(n)})=\eta_v(\sigma)\tau(v_1,\dots,v_n),
\end{align}
alors on peut faire correspondre \`a toute application
$n$-lin\'eaire $deg$-antisym\'etrique $C'$ une application
$n$-lin\'eaire $|~|$-sym\'etrique $C$ en posant $C=\tau C'$ ou
$$
C(v_1,\dots,v_n)=\tau(v_1,\dots,v_n)C'(v_1,\dots,v_n).
$$
Un tel choix a \'et\'e fait dans \cite{[AMM]} o\`u on a pos\'e
$\tau(v_1,\dots, v_n)=\eta_v(\sigma)$ o\`u $\sigma$ est la
permutation rangeant les $v_i$ $deg$-pairs au d\'ebut, sans changer leur ordre 
et les $v_j$ $deg$-impairs en fin sans changer leur ordre.\\

Du fait de la relation que nous venons d'\'etablir entre $[~,~]$ et $Q$, nous 
posons ici
$$
\tau(\alpha_1,\dots,\alpha_n)=(-1)^{\sum_{i=1}^n(n-i)deg(\alpha_i)}.
$$
En fait \ref{(*)} est vrai pour $\sigma=(i,i+1)$ donc pour tout $\sigma$ puisque 
chaque membre de \ref{(*)} d\'efinit une action de $S_n$. On a donc 
$$Q=\tau[~,~].$$\\

\vskip 0.3cm
\begin{defn} {\rm(Cohomologie de Chevalley)}

\

Soit $C'$ une $n$-cocha\^\i ne, c'est \`a dire une application $n$-lin\'eaire 
$deg$-antisym\'etrique de $\left(T_{poly}(\R^d)\right)^n$ dans $T_{poly}(\R^d)$, 
homog\`ene de degr\'e $deg(C')$. Le cobord $\partial'C'$ de $C'$ est par
d\'efinition l'application $n+1$-lin\'eaire $deg$-antisym\'etrique
\begin{align*}
(\partial'&C')(\alpha_0,\dots,\alpha_n)=\cr&\sum_{i=0}^n(-1)^i\varepsilon_{deg
(\alpha)}(i,0\dots\hat{\imath}\dots n)(-1)^{deg(C')deg(\alpha_i)}\left[\alpha_i,
C'(\alpha_0,\dots\widehat{\alpha_i}\dots,\alpha_n)\right]\cr
&-\frac{1}{2}\sum_{i\neq j}\varepsilon_{deg(\alpha)}(i,j,0,\dots\hat{\imath}
\dots\hat{\jmath}\dots n)(-1)^{i+j-1}C'\left([\alpha_i,\alpha_j],\alpha_0,\dots
\widehat{\alpha_i}\dots\widehat{\alpha_j}\dots,\alpha_n\right).
\end{align*}
Une $n$-cocha\^\i ne $C'$ telle que $\partial'C'=0$ est appel\'ee un cocycle, 
une $n$-cocha\^\i ne de la forme $C'=\partial'A'$ est appel\'ee un cobord.
\end{defn}

\

Soit $C'$ une $n$-cocha\^\i ne. Posons $C=\tau C'$ et d\'efinissons $\partial C$
par
$$
\partial C=\tau\partial'C'.$$

\vskip 0.3cm
\begin{lem} {\rm (Cohomologie sym\'etris\'ee)}

\

Une application $n$-lin\'eaire $|~|$-sym\'etrique $C$ est une $n$-cocha\^\i ne, 
son cobord $\partial C$ est donn\'e par
\begin{align*}
(\partial C)&(\alpha_0,\dots,\alpha_n)=\cr&\sum_{i=0}^n\varepsilon_{|\alpha|}
(i,0\dots\hat{\imath}\dots n)(-1)^{|C|(|\alpha_i|-1)}Q\left(\alpha_i,C(\alpha_0,
\dots\widehat{\alpha_i}\dots,\alpha_n)\right)\cr&-\frac{1}{2}\sum_{i\neq j}
\varepsilon_{|\alpha|}(i,j,0,\dots\hat{\imath}\dots\hat{\jmath}\dots
n)C\left(Q(\alpha_i,\alpha_j),\alpha_0,\dots\widehat{\alpha_i}\dots
\widehat{\alpha_j}\dots,\alpha_n\right).
\end{align*}
Ou bien par
\begin{align*}
(\partial C)(\alpha_0,\dots,\alpha_n)&=\sum_{i=0}^n\Big(\varepsilon_{|\alpha|}
(i,0\dots\hat{\imath}\dots n)(-1)^{|C|(|\alpha_i|-1)}\nabla_{\alpha_i}
C(\alpha_0,\dots\widehat{\alpha_i}\dots,\alpha_n)\cr
&+(-1)^{|C|}\varepsilon_{|\alpha|}(0\dots\hat{\imath}\dots n,i)
\nabla_{C(\alpha_0,\dots\widehat{\alpha_i}\dots,\alpha_n)}\alpha_i\Big)\cr
&-\sum_{i\neq j}\varepsilon_{|\alpha|}(i,j,0,\dots\hat{\imath}\dots\hat{\jmath}
\dots n)C\left(\nabla_{\alpha_i}\alpha_j,\alpha_0,\dots\widehat{\alpha_i}\dots
\widehat{\alpha_j}\dots,\alpha_n\right).
\end{align*}
\end{lem}


\

\section{Graphes a\'eriens et cochaines}\label{sec3}

\

Dans cet article, on consid\`ere des graphes de Kontsevich, c'est
\`a dire des graphes $\Gamma$ ayant des sommets a\'eriens
num\'erot\'es $1,\dots,n$, que l'on peut voir comme des points du
demi espace de Poincar\'e $\{\mathfrak{Im}(z)>0\}$ et des sommets
terrestres, num\'erot\'es $\bar{1},\dots,\bar{m}$ que l'on peut
voir comme des points rang\'es sur l'axe r\'eel. De chaque sommet
a\'erien $i$ part $k_i\geq0$ ar\^etes du graphe. Ces ar\^etes sont
des fl\`eches d'extr\'emit\'e soit un sommet a\'erien (on
s'autorise des `petites boucles', c'est \`a dire des ar\^etes de
la forme $\overrightarrow{ii}$) soit un sommet terrestre. Il n'y a
pas d'ar\^ete partant d'un sommet terrestre appel\'e 'pied' du
graphe et il y a exactement une ar\^ete y arrivant, cette ar\^ete
est une jambe du graphe. Il n'y a pas d'ar\^ete multiple (mais on
peut avoir les ar\^etes $\overrightarrow{ij}$ et
$\overrightarrow{ji}$ si $i\neq j$).

Fixons un ordre ${\mathcal O}$ sur les ar\^etes qui soit
compatible avec la num\'erotation des sommets a\'eriens c'est \`a
dire tel que les $k_1$ premi\`eres fl\`eches partent du sommet 1,
les $k_2$ suivantes du sommet 2, etc... On d\'efinit alors une
application $B_{\Gamma,{\mathcal O}}$ $n$-lin\'eaire de
$T_{poly}(\R^d)$ dans lui m\^eme de la fa\c con suivante :

$B_{\Gamma,{\mathcal O}}(\alpha_1,\dots,\alpha_n)$ est nul sauf si
$\alpha_1$ est un $k_1$-tenseur, $\alpha_2$ un
$k_2$-tenseur,\dots,$\alpha_n$ un $k_n$-tenseur.

Dans ce dernier cas, on appelle $Deb(i)$ l'ensemble ordonn\'e
des fl\`eches issues de $i$, elles portent les num\'eros
$a_i=(\sum_{j<i}k_j)+1,a_i+1,\dots, a_i+k_i=\sum_{j\leq i}k_j$. On
pose:
$$
\alpha_i^{Deb(i)}:=\alpha_i^{t_{a_i}t_{a_i+1}\dots t_{a_i+k_i}}.
$$
Pour un sommet a\'erien on note $Fin(i)$ l'ensemble des fl\`eches
arrivant sur $i$. Elles portent les num\'eros $s_1,\dots,s_r$. On
pose:
$$
\partial_{Fin(i)}:=\partial_{t_{s_1}\dots t_{s_r}}.
$$
Pour un pied (un sommet terrestre), on note $Fin(\bar{i})$
l'ensemble des fl\`eches arrivant sur $\bar{i}$. Cet ensemble
contient une seule fl\`eche de num\'ero $s$ et on pose
$$
\partial_{Fin(\bar{i})}:=\partial_{t_s}.
$$
Alors on d\'efinit $B_{\Gamma,{\mathcal O}}(\alpha_1,\dots,\alpha_n)$ par:
$$
B_{\Gamma,{\mathcal O}}(\alpha_1,\dots,\alpha_n)=\sum_{1\leq t_1,\dots, t_{|k|}
\leq d}~~~~\prod_{i=1}^n\partial_{Fin(i)}~~\alpha_i^{Deb(i)}~~
\partial_{Fin(\bar{1})}\wedge\dots\wedge\partial_{Fin(\bar{m})}.
$$
(On a pos\'e $|k|=k_1+\dots+k_n$).

\vskip 0.3cm
\begin{rema}

\

La d\'efinition de l'op\'erateur $B_{\Gamma,{\mathcal O}}$ d\'epend du choix de 
l'ordre compatible ${\mathcal O}$. Changer cet ordre revient \`a multiplier 
$B_{\Gamma,{\mathcal O}}$ par le signe not\'e $\varepsilon({\mathcal O}, 
{\mathcal O}')$ dans \cite{[AGM]} de la permutation de l'ensemble des ar\^etes
faisant passer de ${\mathcal O}$ \`a ${\mathcal O}'$.
\end{rema}

\vskip 0.3cm
On \'etend la d\'efinition de l'op\'erateur $B_{\Gamma,{\mathcal O}}$ aux ordres 
${\mathcal O}'$ non compatibles en posant
$$
B_{\Gamma,{\mathcal O}'}=\varepsilon({\mathcal O},{\mathcal O}')B_{\Gamma, 
{\mathcal O}}.
$$

Si on se restreint aux graphes vectoriels c'est \`a dire aux graphes tels que 
$Deb(i)$ a au plus un \'el\'ement, il y a un seul ordre compatible. On 
num\'erotera alors les fl\`eches par le num\'ero $i$ de leur origine.

Un graphe qui n'a aucun pied et aucune jambe est appel\'e graphe a\'erien. Si 
$\Delta$ est un tel graphe a\'erien, on peut le `compl\'eter' en lui ajoutant 
des pieds et des jambes. On consid\'erera donc tous les graphes $\Gamma$ tels 
que, lorsque l'on retire les pieds et les jambes de $\Gamma$, on retrouve
$\Delta$. Si $(\Delta,{\mathcal O})$ est un graphe orient\'e, on
consid\'erera tous les graphes orient\'es $(\Gamma,{\mathcal O}_\Gamma)$ tels 
que l'ordre induit par ${\mathcal O}_\Gamma$ sur l'ensemble des ar\^etes qui ne 
sont pas des jambes soit exatement ${\mathcal O}$. On notera cette propri\'et\'e 
$(\Gamma, {\mathcal O}_\Gamma)\supset(\Delta,{\mathcal O})$.

Plus pr\'ecis\'ement, on reprendra la d\'efinition de \cite{[AGM]}. Partant de 
$(\Delta,{\mathcal O})$ o\`u ${\mathcal O}$ est compatible et de $m$ pieds 
num\'erot\'es $\bar{1},\dots,\bar{m}$, on peut construire un graphe $\Gamma$ en 
ajoutant $m$ jambes \`a $\Delta$ (une pour chaque pied). On peut d\'efinir un 
ordre ${\mathcal O}_0$ sur les ar\^etes de $\Gamma$ en rangeant les jambes 
apr\`es les ar\^etes a\'eriennes, dans l'ordre des pieds. Si $k_1,\dots,k_n$ 
sont les nombres d'ar\^etes de $\Gamma$ issues des sommets $1,\dots,n$ et 
$\ell_1,\dots,\ell_n$ le nombre d'ar\^etes de $\Delta$ issues des sommets $1,
\dots,n$, il y a $\displaystyle{\frac{k!}{\ell!}}$ ordres compatibles possibles
${\mathcal O}_\Gamma$ tels que $(\Gamma, {\mathcal O}_\Gamma)\supset(\Delta,
{\mathcal O})$. On posera:
$$
C_{\Delta,{\mathcal O}}=\sum_{m=0}^\infty~\frac{1}{m!}~
\sum_{\begin{smallmatrix}(\Gamma,{\mathcal O}_\Gamma)\supset(\Delta,{\mathcal
O})\cr\#\{pieds~de~\Gamma\}=m\end{smallmatrix}}\frac{\ell!}{k!}
\varepsilon({\mathcal O}_\Gamma,{\mathcal O}_0)B_{\Gamma,{\mathcal O}_\Gamma}.
$$
Lorsque l'on se restreint aux graphes vectoriels, cette somme est finie. La 
somme ci-dessus ne contient que des graphes $\Gamma$ ayant $m$ pieds avec:
$$
0\leq m\leq n-\vert\ell\vert.
$$
Comme il n'y a qu'un ordre compatible sur les ar\^etes de $\Delta$.\\

Les cocha\^\i nes $C_\delta$ \'etudi\'ees dans ce papier sont les 
sym\'etris\'ees des applications $C_{\Delta,{\mathcal O}}$, c'est \`a dire des
op\'erateurs associ\'es \`a des combinaisons lin\'eaires sym\'etriques de
graphes
$$
\delta=\sum_{\Delta,{\mathcal O}}a_{\Delta,{\mathcal O}}(\Delta,{\mathcal O}).
$$
Si $\sigma$ est une premutation de $n$ \'el\'ements, $\sigma$ agit sur $\Delta$ 
en permutant ses sommets et ses ar\^etes par paquets, en gardant l'ordre des 
ar\^etes issues d'un m\^eme sommet. On sait alors (\cite{[AGM]}, 
Proposition 4.4) que:
$$
C_{\Delta,{\mathcal O}}(\alpha_{\sigma(1)},\dots,\alpha_{\sigma(n)})=
\varepsilon_{|Deb(\Delta)|}(\sigma)C_{\Delta,{\mathcal O}}(\alpha_1,\dots,
\alpha_n).
$$
On dira donc que $\delta$ est sym\'etrique si $C_\delta=\sum_{\Delta,{\mathcal 
O}}a_{\Delta,{\mathcal O}}C_{\Delta,{\mathcal O}}$ l'est, c'est \`a dire si, 
pour tout $(\Delta,{\mathcal O})$,
$$
a_{\sigma(\Delta),\sigma({\mathcal O})}=\varepsilon_{|Deb(\Delta)|}
a_{\Delta,{\mathcal O}}.
$$

Par exemple la roue simple $\Delta$ de longueur 3 est un graphe a\'erien \`a 3 
sommets num\'erot\'es 1, 2, 3 et ayant pour ar\^etes
$\{\overrightarrow{12},\overrightarrow{23},\overrightarrow{31}\}$. C'est un 
graphe vectoriel. La sym\'etrisation $\delta$ de $(\Delta,{\mathcal O})$ est:
$$
\delta=\sum_{\sigma\in S_3}\varepsilon(\sigma)(\sigma(\Delta),\sigma({\mathcal
O}))=3\times\{\overrightarrow{12},\overrightarrow{23},\overrightarrow{31}\}
-3\times\{\overrightarrow{13},\overrightarrow{32},\overrightarrow{21}\}.
$$
Cette sym\'etrisation d\'efinit une application $C_\delta$ qui envoie trois 
champs de vecteurs $\alpha_1$, $\alpha_2$, $\alpha_3$ sur la fonction
$$
C_\delta(\alpha_1,\alpha_2,\alpha_3)=3\partial_{i_3}\alpha_1^{i_1}\partial_{i_1}
\alpha_2^{i_2}\partial_{i_2}\alpha_3^{i_3}-3\partial_{i_2}\alpha_1^{i_1}
\partial_{i_1}\alpha_3^{i_3}\partial_{i_3}\alpha_2^{i_2}.
$$


\

\

\section{L'op\'erateur de cobord sur les graphes vectoriels}\label{sec4}

\

Dans ce paragraphe, nous d\'efinissons directement sur les combinaisons
lin\'eaires $\delta$ sym\'etriques de graphes a\'eriens vectoriels un 
op\'erateur de cobord $\partial$ correspondant \`a l'op\'erateur de cobord 
pour les cocha\^\i nes sym\'etriques $C_\delta$ d\'efinies par $\delta$. Dans
\cite{[AGM]}, il est montr\'e qu'un tel op\'erateur d\'efini sur les graphes 
existe. Dans le cas des graphes vectoriels, son expression peut \^etre 
simplifi\'ee.

\vskip 0.3cm
\begin{defn} {\rm(Eclatement propre d'un sommet)}

\

Soit $\Delta$ un graphe vectoriel a\'erien de sommets num\'erot\'es $1,\dots,n$.
Fixons $i$ et $j$ tels que $0\leq i,j\leq n$. On construit une famille de 
graphes $\Delta'_{ji}$ (resp. $\Delta'_{ij}$) de la fa\c con suivante:

On renum\'erote les sommets de $\Delta$ en $0,\dots\hat{\jmath}\dots,n$ (en 
gardant leur ordre initial),\\

\noindent
\underbar{Si le sommet $i$ est un point isol\'e ($Deb(i)=Fin(i)=\emptyset$)}

On ajoute le sommet $j$ et la fl\`eche $\overrightarrow{ji}$ (resp. 
$\overrightarrow{ij}$).\\

\noindent
\underbar{Si le sommet $i$ n'est pas un point isol\'e}

On ajoute un sommet $j$ et une fl\`eche $\overrightarrow{ji}$ (resp. 
$\overrightarrow{ij}$),

Si une fl\`eche $\overrightarrow{ia}$ partait de $i$ dans $\Delta$, on la garde 
dans $\Delta'_{ji}$ (resp. on la remplace par $\overrightarrow{ja}$ dans 
$\Delta'_{ij}$),

On r\'epartit les fl\`eches du graphe $\Delta$ arrivant sur le sommet $i$ entre 
les sommets $j$ et $i$ du nouveau graphe $\Delta'_{ji}$ (resp. $\Delta'_{ij}$)
proprement, c'est \`a dire de telle fa\c con que:
$$\inf\left(|Deb(i)|+|Fin(i)|,|Deb(j)|+|Fin(j)|\right)>1.$$
\end{defn}

\vskip 0.3cm
Remarquons que:\\

Si $i$ est un sommet isol\'e, il y a un seul graphe 
$\Delta'_{ji}$.

Si $i$ est tel que $|Deb(i)|+|Fin(i)|=1$ alors il n'y a aucun graphe 
$\Delta'_{ji}$.

Si $i$ est un sommet tel que $|Deb(i)|=|Fin(i)|=1$, il y a un seul 
$\Delta'_{ji}$.

En g\'en\'eral si $|Deb(i)|=1$, il y a $2^{|Fin(i)|}-1$ graphes $\Delta'_{ji}$,
si $|Deb(i)|=0$ et $|Fin(i)|\geq1$, il y a $2^{|Fin(i)|}-2$ graphes 
$\Delta'_{ji}$.\\

Le m\^eme r\'esultat est valable pour $\Delta'_{ij}$.\\

Pour repr\'esenter la famille des graphes qu'on vient de d\'efinir, on notera:
$$\Delta'_{ji}\rightarrow^{prop}_i\Delta\quad (\hbox{resp.}\Delta'_{ij}
\rightarrow^{prop}_i\Delta).$$\\

\vskip 0.3cm
\begin{prop}\label{partial} {\rm (L'op\'erateur $\partial$ sur les graphes)}

\

Soit $\delta=\sum_\Delta a_\Delta\Delta$ une combinaison lin\'eaire sym\'etrique 
homog\`ene de graphes a\'eriens vectoriels. Alors
$$
\partial C_\delta=\sum_\Delta a_\Delta C_{\partial\Delta}
$$
avec
$$
\partial\Delta=-\sum_{j\neq i}\varepsilon_{\{|Deb(0)|,\dots,|Deb(n)|\}}(j,0,
\dots\hat{\jmath}\dots,n)\sum_{\Delta'_{ji}\mathop{\rightarrow^{prop}_i}\Delta}
\Delta'_{ji}.
$$
On peut \'ecrire $\partial\Delta$ autrement. Notons $\varepsilon_{|Deb|}$ la 
quantit\'e $\varepsilon_{\{|Deb(0)|,\dots,|Deb(n)|\}}$, alors
\begin{align*}
\partial\Delta=-&\sum_{i<j}\Big[\varepsilon_{|Deb|}(j,0,\dots\hat{\jmath}\dots,
n)\sum_{\Delta'_{ji}\rightarrow^{prop}_i\Delta}\Delta'_{ji}\cr
&+\big[\left(1-|Deb(j)|\right)\varepsilon_{|Deb|}(i,0,\dots\hat{\imath}\dots,n)-
|Deb(j)|\varepsilon_{|Deb|}(j,0,\dots\hat{\jmath}\dots,n)\big]\cr
&\times\sum_{\Delta'_{ij}\rightarrow^{prop}_i\Delta}\Delta'_{ij}\Big]
\end{align*}

\end{prop}

\vskip 0.3cm
\noindent
{\bf Preuve}

On sait que:
\begin{align*}
(\partial C_\delta)(\alpha_0,\dots,\alpha_n)&=\sum_{j=0}^n
\varepsilon_{|\alpha|}(j,0\dots\hat{\jmath}\dots n)(-1)^{|C|(|\alpha_j|-1)}
\nabla_{\alpha_j}C_\delta(\alpha_0,\dots\widehat{\alpha_j}\dots,\alpha_n)\cr
&+\sum_{i=0}^n(-1)^{|C_\delta|}\varepsilon_{|\alpha|}(0\dots\hat{\imath}\dots 
n,i)\nabla_{C_\delta(\alpha_0,\dots\widehat{\alpha_i}\dots,\alpha_n)}\alpha_i\cr
&-\sum_{i\neq j}\varepsilon_{|\alpha|}(j,i,0,\dots\hat{\imath}\dots\hat{\jmath}
\dots n)C_\delta\left(\nabla_{\alpha_j}\alpha_i,\alpha_0,\dots\widehat{\alpha_i}
\dots\widehat{\alpha_j}\dots,\alpha_n\right)\cr
&=(I)+(II)-(III).
\end{align*}
Puisque la cohomologie de $C_\delta$ est d\'etermin\'ee par sa composante dans 
les fonctions (\cite{[AGM]}, Proposition 5.2), on ne regarde que les termes qui
sont des op\'erateurs 0 diff\'erentiels (des fonctions) dans cette expression.

Les termes de la premi\`ere somme n'apparaissent que si $\alpha_j$ est un champ 
de vecteurs $\alpha_j=\sum_\ell\alpha_j^\ell\partial_\ell$. On a
$$
(I)=\sum_{j\neq i}\varepsilon_{|\alpha|}(j,0,\dots\hat{\jmath}\dots,n)
\sum_\ell\alpha_j^\ell C_\delta(\alpha_0,\dots\widehat{\alpha_j}\dots,
\partial_\ell\alpha_i,\dots,\alpha_n)
$$
Dans le second terme, $|C_\delta|$ est le nombre $|\delta|$ de fl\`eches des
graphes de $\delta$. En revenant \`a la d\'efinition de $C_\delta$, ces termes 
apparaissent lorsque $C_\delta(\alpha_0,\dots\widehat{\alpha_i}\dots,\alpha_n)$
est un champ de vecteurs. Pour chaque $\Delta$ de $\delta$, on choisit un sommet 
$j$ tel que $Deb(\alpha_j)=0$. Alors:
\begin{align*}
C_\Delta&(\alpha_0,\dots\widehat{\alpha_i}\dots,\alpha_n)\cr
&=\sum_{(\Gamma,{\mathcal O}_\Gamma)\supset(\Delta,{\mathcal O})}\varepsilon(
{\mathcal O}_\Gamma,{\mathcal O}_0)\sum_{j,\ell}B_{(\Gamma,{\mathcal O})}(
\alpha_0,\dots,\alpha_j^\ell,\dots\widehat{\alpha_i}\dots,\alpha_n)\partial_\ell
\cr
&=\sum_{(\Gamma,{\mathcal O}_\Gamma)\supset(\Delta,{\mathcal O})}
\varepsilon_{|\alpha|}(0,\dots\hat{\imath}\dots\hat{\jmath}\dots,n,j)
\sum_{j,\ell}B_{(\Gamma,{\mathcal O})}(\alpha_0,\dots,\alpha_j^\ell,\dots
\widehat{\alpha_i}\dots,\alpha_n)\partial_\ell.
\end{align*}
Donc on obtient puisqu'ici $|\alpha_i|=0$,
\begin{align*}
(II)&=\sum_{j\neq i}(-1)^{|\delta|}\varepsilon_{|\alpha|}(0,\dots\hat{\imath}
\dots\hat{\jmath}\dots,n,j)\sum_\ell C_\delta(\alpha_0,\dots,\alpha_j^\ell,\dots
\widehat{\alpha_i}\dots,\alpha_n)\partial_\ell\alpha_i\cr
&=\sum_{j\neq i}\varepsilon_{|\alpha|}(j,0,\dots\hat{\jmath}
\dots,n)\sum_\ell C_\delta(\alpha_0,\dots,\alpha_j^\ell,\dots
\widehat{\alpha_i}\dots,\alpha_n)\partial_\ell\alpha_i\cr
&=\sum_{j\neq i}\varepsilon_{|\alpha|}(j,0,\dots\hat{\jmath}
\dots,n)\sum_\ell C_\delta(\alpha_0,\dots,\underbrace{\alpha_j^\ell}_{(i)},\dots
\widehat{\alpha_j}\dots,\alpha_n)\partial_\ell\alpha_i.
\end{align*}
Enfin la derni\`ere somme s'\'ecrit:
\begin{align*}
(III)&=\sum_{i\neq j}\varepsilon_{|\alpha|}(j,i,0,\dots\hat{\imath}\dots\hat
{\jmath}\dots,n)C_\delta\left(\nabla_{\alpha_j}\alpha_i,\alpha_0,\dots
\widehat{\alpha_i}\dots\widehat{\alpha_j}\dots,\alpha_n\right)\cr
&=\sum_{i\neq j}\varepsilon_{|\alpha|}(j,0,\dots\hat{\jmath}\dots,n)\sum_\ell 
C_\delta(\alpha_0,\dots,\underbrace{\alpha_j^\ell\partial_\ell\alpha_i
}_{(i)},\dots\widehat{\alpha_j}\dots,\alpha_n).
\end{align*}
Appliquons la r\`egle de Leibniz pour les d\'erivations multiples dans le terme 
$(III)$, on obtient pour chaque $\Delta$ de $\delta$:
\begin{align*}
C_\Delta&(\alpha_0,\dots,\underbrace{\alpha_j^\ell\partial_\ell\alpha_i
}_{(i)},\dots\widehat{\alpha_j}\dots,\alpha_n)\cr
&=\prod_{t\neq i}\partial_{Fin(t)}\alpha_t^{Deb(t)}.\partial_{Fin(i)}\left(
\alpha_j^\ell\partial_\ell\alpha_i^{Deb(i)}\right)\cr
&=\prod_{t\neq i}\partial_{Fin(t)}\alpha_t^{Deb(t)}.\left(\sum_{A\subset Fin(i)}
\partial_A\alpha_j^\ell\partial_{\{\ell\}\cup Fin(i)\setminus A}
\alpha_i^{Deb(i)}
\right).
\end{align*}
La somme $(I)$ correspond exactement au cas $A=\emptyset$, la somme $(II)$ au 
cas $A=Fin(i)$. Il y a simplification. Posons $\varepsilon_{|Deb|}=
\varepsilon_{\{|Deb(0)|,\dots,|Deb(n)|\}}$ et
$$
\partial C_\Delta=-\sum_{j\neq i}\varepsilon_{|Deb|}(j,0,\dots
\hat{\jmath}\dots,n)\partial_{ji}C_\Delta=-\sum_{j\neq i}\varepsilon_{|Deb|}
(j,0,\dots\hat{\jmath}\dots,n)C_{\partial_{ji}\Delta},
$$
on peut \'ecrire:\\

\noindent
\underbar{Si $Deb(i)=Fin(i)=\emptyset$}
$$
\partial_{ji}\Delta=\Delta'_{ji}=\sum_{\Delta'_{ji}\rightarrow^{prop}_i\Delta}
\Delta'_{ji}.
$$\\

\noindent
\underbar{Si $Deb(i)=\emptyset$ et $|Fin(i)|>1$}
\begin{align*}
\partial_{ji}C_\Delta&=\prod_{t\neq i}\partial_{Fin(t)}\alpha_t^{Deb(t)}.\left(
\sum_{\begin{smallmatrix}A\subset Fin(i)\cr A\neq\emptyset\cr A\neq Fin(i)
\end{smallmatrix}}\partial_A\alpha_j^\ell\partial_{\{\ell\}\cup Fin(i)\setminus 
A}\alpha_i\right)\cr
&=\sum_{\Delta'_{ji}\rightarrow^{prop}_i\Delta}C_{\Delta'_{ji}}.
\end{align*}\\

\noindent
\underbar{Si $Deb(i)\neq\emptyset$ et $Fin(i)\neq\emptyset$}
\begin{align*}
\partial_{ji}C_\Delta&=\prod_{t\neq i}\partial_{Fin(t)}\alpha_t^{Deb(t)}.\left(
\sum_{\begin{smallmatrix}A\subset Fin(i)\cr A\neq\emptyset
\end{smallmatrix}}\partial_A\alpha_j^\ell\partial_{\{\ell\}\cup Fin(i)\setminus 
A}\alpha_i^{Deb(i)}\right)\cr
&=\sum_{\Delta'_{ji}\rightarrow^{prop}_i\Delta}C_{\Delta'_{ji}}.
\end{align*}\\

\noindent
\underbar{Si $|Deb(i)|+|Fin(i)|=1$}
$$
\partial_{ji}\Delta=0=\sum_{\Delta'_{ji}\rightarrow^{prop}_i\Delta}
\Delta'_{ji}.
$$\\
Finalement on a donc bien:
$$
\partial\Delta=-\sum_{j\neq i}\varepsilon_{|Deb|}(j,0,\dots\hat{\jmath}\dots,n)
\sum_{\Delta'_{ji}\rightarrow^{prop}_i\Delta}\Delta'_{ji}.
$$\\

D'autre part, on peut regrouper les termes de $(III)$ autrement:
\begin{align*}
(III)&=\sum_{i<j}\varepsilon_{|\alpha|}(j,i,0,\dots\hat{\imath}\dots\hat
{\jmath}\dots,n)C_\delta\left(\nabla_{\alpha_j}\alpha_i,\alpha_0,\dots
\widehat{\alpha_i}\dots\widehat{\alpha_j}\dots,\alpha_n\right)\cr
&\hskip 1cm+\sum_{i<j}
\varepsilon_{|\alpha|}(i,j,0,\dots\hat{\imath}\dots\hat{\jmath}\dots,n)C_\delta
\left(\nabla_{\alpha_i}\alpha_j,\alpha_0,\dots\widehat{\alpha_i}\dots
\widehat{\alpha_j}\dots,\alpha_n\right)\cr
\end{align*}
ou
\begin{align*}
&(III)=\sum_{i<j}\varepsilon_{|\alpha|}(j,0,\dots\hat{\jmath}\dots,n)\sum_\ell 
C_\delta(\alpha_0,\dots,\underbrace{\alpha_j^\ell\partial_\ell\alpha_i
}_{(i)},\dots\widehat{\alpha_j}\dots,\alpha_n)\cr
&\hskip 0.5cm+\sum_{i<j}(-1)^{|\alpha_i||\alpha_j|}\varepsilon_{|\alpha|}(j,i,0,
\dots\hat{\imath}\dots\hat{\jmath}\dots,n)\sum_\ell C_\delta
\left(\alpha_i^\ell\partial_\ell\alpha_j,\alpha_0,\dots\widehat{\alpha_i}\dots
\widehat{\alpha_j}\dots,\alpha_n\right).
\end{align*}
On transforme le second terme ainsi:\\

\noindent
Si $|\alpha_j|=1$, ce terme est:
$$
-\varepsilon_{|\alpha|}(j,0,\dots\hat{\jmath}\dots,n)\sum_\ell 
C_\delta(\alpha_0,\dots,\underbrace{\alpha_i^\ell\partial_\ell\alpha_j}_{(i)},
\dots\widehat{\alpha_j}\dots,\alpha_n).
$$\\

\noindent
Si $|\alpha_j|=0$, ce terme est:
$$
\varepsilon_{|\alpha|}(i,0,\dots\hat{\imath}\dots,n)\sum_\ell 
C_\delta(\alpha_0,\dots,\underbrace{\alpha_i^\ell\partial_\ell\alpha_j}_{(i)},
\dots\widehat{\alpha_j}\dots,\alpha_n).
$$
En regroupant les deux cas, on obtient le facteur annonc\'e:
\begin{align*}
\big[&\left(1-|\alpha_j|\right)\varepsilon_{|\alpha|}(i,0,\dots\hat{\imath}
\dots,n)-|\alpha_j|\varepsilon_{|\alpha|}(j,0,\dots\hat{\jmath}\dots,n)\big]\cr
&\hskip 2cm C_\delta(\alpha_0,\dots,\underbrace{\alpha_i^\ell\partial_\ell
\alpha_j}_{(i)},\dots\widehat{\alpha_j}\dots,\alpha_n).
\end{align*}

De m\^eme pour $(I)$ et $(II)$, on obtient
\begin{align*}
(I)&=\sum_{i<j}\varepsilon_{|\alpha|}(j,0,\dots\hat{\jmath}\dots,n)
\sum_\ell\alpha_j^\ell C_\delta(\alpha_0,\dots,\underbrace{\partial_\ell
\alpha_i}_{(i)},\dots\widehat{\alpha_j}\dots,\alpha_n)\cr
&\hskip 1cm+\big[\left(1-|\alpha_j|\right)\varepsilon_{|\alpha|}(i,0,\dots
\hat{\imath}\dots,n)-|\alpha_j|\varepsilon_{|\alpha|}(j,0,\dots\hat{\jmath}
\dots,n)\big]\cr
&\hskip 1cm\sum_\ell\alpha_i^\ell C_\delta(\alpha_0,\dots,\underbrace{
\partial_\ell\alpha_j}_{(i)},\dots\widehat{\alpha_j}\dots,\alpha_n)
\end{align*}
et
\begin{align*}
(II)&=\sum_{i<j}\varepsilon_{|\alpha|}(j,0,\dots\hat{\jmath}\dots,n)
\sum_\ell C_\delta(\alpha_0,\dots,\underbrace{\alpha_j^\ell}_{(i)},
\dots\widehat{\alpha_j}\dots,\alpha_n)\partial_\ell\alpha_i\cr
&\hskip 1cm+\varepsilon_{|\alpha|}(i,0,\dots\hat{\imath}\dots,n)\sum_\ell 
C_\delta(\alpha_0,\dots,\underbrace{\alpha_i^\ell}_{(i)},
\dots\widehat{\alpha_j}\dots,\alpha_n)\partial_\ell\alpha_j.
\end{align*}

En faisant le m\^eme calcul que ci-dessus, on simplifie les termes qui se 
correspondent dans $(I)$, $(II)$ et $(III)$ et on obtient la formule annonc\'ee.

\vskip 0.3cm
\begin{rema}

\

On a retrouv\'e le cobord d\'efini dans \cite{[AGM]}. En particulier, si 
$\delta$ est sym\'etrique, $\partial\delta$ est aussi sym\'etrique.
\end{rema}

On pose bien entendu

\vskip 0.3cm
\begin{defn} {\rm (Espaces de cohomologie)}

\

Une combinaison lin\'eaire sym\'etrique $\delta$ de graphes vectoriels a\'eriens 
est un cocycle si $\partial\delta=0$, un cobord s'il existe une combinaison
lin\'eaire sym\'etrique $\beta$ telle que $\delta=\partial\beta$.

L'espace $Z^n$ des $n$-cocycles est l'espace des combinaisons sym\'etriques de 
graphes vectoriels ayant $n$ sommets et qui sont des cocycles.

L'espace $B^n$ des $n$-cobords est l'espace des combinaisons de graphes qui sont 
les cobords de combinaisons sym\'etriques de graphes ayant $n-1$ sommets.

Le $n^{\grave{e}me}$ espace de cohomologie des graphes $H^n$ est le quotient de
l'espace $Z^n$ par l'espace $B^n$.
\end{defn}


\

\section{Symbole d'un graphe}

\

Soit $\Delta$ un graphe a\'erien vectoriel de sommets num\'erot\'es 
$(1,\dots,n)$. On peut distinguer six classes de sommets $i$:\\

Classe 1: les sommets $i$ tels que $|Fin(i)|>1$ et $|Deb(i)|=0$. On appellera 
ordre de $i$ le symbole $r_i$ o\`u $r_i=|Fin(i)|$,\\

Classe 2: les sommets $i$ tels que $|Fin(i)|>1$ et $|Deb(i)|=1$. On appellera 
ordre de $i$ le symbole $r_i^+$ o\`u $r_i=|Fin(i)|$,\\

Classe 3: les sommets $i$ tels que $|Fin(i)|=1$ et $|Deb(i)|=1$. On appellera 
ordre de $i$ le symbole $1^+$,\\

Classe 4: les sommets $i$ tels que $|Fin(i)|=1$ et $|Deb(i)|=0$. On appellera 
ordre de $i$ le symbole $1$,\\

Classe 5: les sommets $i$ tels que $|Fin(i)|=0$ et $|Deb(i)|=0$. On appellera 
ordre de $i$ le symbole $0$,\\

Classe 6: les sommets $i$ tels que $|Fin(i)|=0$ et $|Deb(i)|=1$. On appellera 
ordre de $i$ le symbole $0^-$.\\

On ordonne les ordres des sommets en posant:
$$
r_i>r_j^+>1^+>1>0>0^-~~\hbox{et}~~r_i^+\geq r_{i'}^+\Leftrightarrow r_i\geq 
r_{i'}.\eqno{(*)}
$$
L'ordre ${\mathcal O}(\Delta)$ d'un graphe $\Delta$ est le mot form\'e par les
ordres de ses sommets:
$$
{\mathcal O}(\Delta)=\left({\mathcal O}(1),\dots,{\mathcal O}(n)\right).
$$
On ordonne les ordres des graphes en utilisant l'ordre lexicographique, en 
respectant $(*)$.

Si $\delta=\sum_\Delta a_\Delta\Delta$ est une combinaison lin\'eaire 
sym\'etrique de graphes vectoriels, on d\'efinit l'ordre de $\delta$ par:
$$
{\mathcal O}(\delta)=Max\{{\mathcal O}(\Delta),~~a_\Delta\neq0\}
$$
et on appellera symbole de $\delta$ la combinaison lin\'eaire non sym\'etrique:
$$
\sigma_\delta=\sum_{\begin{smallmatrix}\Delta\cr{\mathcal O}(\Delta)={\mathcal 
O}(\delta)\end{smallmatrix}}a_\Delta\Delta.
$$
Comme $\delta$ est sym\'etrique, son ordre a la forme:
$$
{\mathcal O}(\delta)=(r_1,\dots,r_{k_0-1},r_{k_0}^+,\dots,r_{k_1-1}^+,
\underbrace{1^+}_{(k_1)},\dots,1^+,\underbrace{1}_{(k_2)},\dots,1,
\underbrace{0}_{(k_3)},\dots,0,\underbrace{0^-}_{(k_4)},\dots,0^-)
$$
avec
$$
r_1\geq r_2\dots\geq r_{k_0-1},~~r_{k_0}^+\geq\dots\geq r_{k_1-1}^+
$$
(on peut avoir $k_0=1$ ou $k_1=k_0$, etc\dots)\\

\vskip 0.3cm
\begin{prop} {\rm (Le symbole de $\partial\delta$)}

\

Soit $\delta$ une combinaison lin\'eaire sym\'etrique de graphes a\'eriens 
vectoriels, d'ordre
\begin{align*}
{\mathcal O}&(\delta)=\cr&(r_1,\dots,r_{k_0-1},r_{k_0}^+,\dots,r_{k_1-1}^+,
\underbrace{1^+}_{(k_1)},\dots,1^+,\underbrace{1}_{(k_2)},\dots,1,
\underbrace{0}_{(k_3)},\dots,0,\underbrace{0^-}_{(k_4)},\dots,0^-).
\end{align*}
Alors chaque graphe $\Delta'$ apparaissant dans $\partial\delta$ est d'ordre
au plus:
\begin{align*}
&{\mathcal O}(\delta)\oplus1^+=\cr
&(r_1,\dots,r_{k_0-1},r_{k_0}^+,\dots,r_{k_1-1}^+,
\underbrace{1^+}_{(k_1)},\dots,1^+,\underbrace{1}_{(k_2+1)},\dots,1,
\underbrace{0}_{(k_3+1)},\dots,0,\underbrace{0^-}_{(k_4+1)},\dots,0^-).
\end{align*}\\

Si ${\mathcal O}(\partial\delta)={\mathcal O}(\delta)\oplus1^+$, alors le 
symbole de $\partial\delta$ est:
\begin{align*}
\sigma_{\partial\delta}&=-\sum_{\Delta\in\sigma_\delta}a_\Delta
\sum_{\begin{smallmatrix}i<j\cr0\leq i<k_0\cr k_1\leq j\leq k_2
\end{smallmatrix}}\varepsilon_{|Deb|}(j,0,\dots\hat{\jmath}\dots,n)
\sum_{\Delta'_{ji}\rightarrow^{prop}_i\Delta}
\Delta'_{ji}\cr
&\hskip 3cm+\sum_{\begin{smallmatrix}i<j\cr
k_0\leq i<k_1\cr k_1\leq j\leq k_2\end{smallmatrix}}\varepsilon_{|Deb|}(j,0,
\dots\hat{\jmath}\dots,n)\Big[
\sum_{\Delta'_{ji}\rightarrow^{prop}_i\Delta}\Delta'_{ji}-\sum_{\Delta'_{ij}
\rightarrow^{prop}_i\Delta}\Delta'_{ij}\Big]\cr
&\hskip 3cm+\sum_{\begin{smallmatrix}i<j\cr k_1\leq i<k_2\cr k_1\leq j\leq k_2
\end{smallmatrix}}\varepsilon_{|Deb|}(j,0,\dots\hat{\jmath}\dots,n)\Big[
\sum_{\Delta'_{ji}\rightarrow^{prop}_i\Delta}\Delta'_{ji}-\sum_{\Delta'_{ij}
\rightarrow^{prop}_i\Delta}\Delta'_{ij}\Big].
\end{align*}
\end{prop}

\vskip 0.3cm
\noindent
{\bf Preuve}

Si $\Delta$ est un graphe de $\delta$ qui n'appara\^\i t pas dans le symbole 
$\sigma_\delta$, alors les ordres des graphes $\Delta'_{ij}$ et $\Delta'_{ji}$
qui se contractent sur $\Delta$ sont tous strictement plus petits que ${\mathcal 
O}(\delta)\oplus1^+$. Regardons donc seulement les $\Delta$ de $\sigma_\delta$.

Fixons un couple $(i,j)$ avec $i<j$. Il est clair que dans la d\'ecomposition de 
$\partial\Delta$, les graphes $\Delta'=\Delta'_{ji}$ ou $\Delta'=\Delta'_{ij}$ 
sont d'ordre:

Si $0\leq i<k_0$, alors
$$
{\mathcal O}(\Delta'_{ji})=(r_1,\dots,\underbrace{(r_i-k+1)}_{(i)},\dots,
\underbrace{k^+}_{(j)},\dots)\qquad (1\leq k<r_i).
$$
Donc ${\mathcal O}(\Delta'_{ji})\leq{\mathcal O}(\delta)\oplus1^+$ et 
l'\'egalit\'e n'est vraie que si $k=1$ et $k_1\leq j\leq k_2$. Il y a 
$r_i$ graphes $\Delta'_{ji}$ dans ce cas. D'autre part:
$$
{\mathcal O}(\Delta'_{ij})=(r_1,\dots,\underbrace{k^+}_{(i)},\dots,
\underbrace{(r_i-k+1)}_{(j)},\dots)\qquad (1\leq k<r_i).
$$
Donc ${\mathcal O}(\Delta'_{ij})<{\mathcal O}(\delta)\oplus1^+$. Il n'y a 
aucun graphe $\Delta'_{ij}$ dans ce cas.

Si $k_0\leq i<k_1$, alors
$$
{\mathcal O}(\Delta'_{ji})=(r_1,\dots,\underbrace{(r_i-k+1)^+}_{(i)},\dots,
\underbrace{k^+}_{(j)},\dots) \qquad(1\leq k\leq r_i).
$$
Donc ${\mathcal O}(\Delta'_{ji})\leq{\mathcal O}(\delta)\oplus1^+$ et 
l'\'egalit\'e n'est vraie que si  $k=1$ et $k_1\leq j\leq k_2$. Il y a $r_i$
graphes $\Delta'_{ji}$ dans ce cas. D'autre part:
$$
{\mathcal O}(\Delta'_{ij})=(r_1,\dots,\underbrace{k^+}_{(i)},\dots,
\underbrace{(r_i-k+1)^+}_{(j)},\dots) \qquad(1\leq k\leq r_i).
$$
Donc ${\mathcal O}(\Delta'_{ij})\leq{\mathcal O}(\delta)\oplus1^+$ et 
l'\'egalit\'e n'est vraie que si $k=r_i$ et $k_1\leq j\leq k_2$. Il y a un seul
graphe $\Delta'_{ij}$ dans ce cas.

Si $k_1\leq i<k_2$, alors
$$
{\mathcal O}(\Delta')=(r_1,\dots,\underbrace{1^+}_{(i)},\dots,
\underbrace{1^+}_{(j)},\dots).$$
Donc pour et seulement pour $i<j\leq k_2$, ${\mathcal O}(\Delta')={\mathcal O}
(\delta)\oplus1^+$ et il y a un seul graphe $\Delta'_{ji}$ et un seul graphe 
$\Delta'_{ij}$ dans ce cas.

Si $k_2<i$, pour tout $j>i$, on a ${\mathcal O}(\Delta')<{\mathcal O}
(\delta)\oplus1^+$.\\

D\'efinissons maintenant l'op\'erateur d'homotopie.

\vskip 0.3cm
\begin{defn} {\rm (L'homotopie)}

\

Soit $\Delta$ un graphe vectoriel de sommets $(0,\dots,n)$, ayant des sommets 
$i$ tels que $|Deb(i)|=|Fin(i)|=1$ (des sommets $i$ d'ordre $1^+$). On d\'efinit 
le graphe $h(\Delta)$ ainsi:

On consid\`ere $i_0$, le plus grand des indices $i$ d'ordre $1^+$. La fl\`eche 
issue de $i_0$ est $\overrightarrow{i_0a}$. Le graphe $h(\Delta)$ est le graphe 
de sommets $(0,\dots\hat{\imath_0}\dots,n)$ obtenu en comprimant la fl\`eche
$\overrightarrow{i_0a}$ et en identifiant les sommets $i_0$ et $a$ au sommet $a
$.

Si $\Delta$ n'a pas de sommets d'ordre $1^+$, on pose $h(\Delta)=0$.

On prolonge $h$ lin\'eairement \`a l'espace des combinaisons lin\'eaires de 
graphes.
\end{defn}

\vskip 0.3cm
\begin{prop} {\rm (Symbole et homotopie)}

\

Soit $\delta$ une combinaison lin\'eaire sym\'etrique de graphes vectoriels 
d'ordre
\begin{align*}
{\mathcal O}&(\delta)=\cr&(r_1,\dots,r_{k_0-1},r_{k_0}^+,\dots,r_{k_1-1}^+,
\underbrace{1^+}_{(k_1)},\dots,1^+,\underbrace{1}_{(k_2)},\dots,1,
\underbrace{0}_{(k_3)},\dots,0,\underbrace{0^-}_{(k_4)},\dots,0^-).
\end{align*}
Par abus de notations, on pose:
$$
\sigma_{\partial\delta}=\sum_{\begin{smallmatrix}\Delta'\in\partial\delta\cr
{\mathcal O}(\Delta')={\mathcal O}(\delta)\oplus1^+\end{smallmatrix}}a_{\Delta'}
\Delta'
$$
si $\partial\delta=\sum_{\Delta'}a_{\Delta'}\Delta'$.

Notons enfin $\delta(0,\dots\hat{k_2}\dots,n)$ la combinaison lin\'eaire de
graphes $\Delta$ dont on a renum\'erot\'e les sommets en $(0,\dots\hat{k_2}\dots,
n)$. Alors
\begin{align*}
h(\sigma_{\partial\delta})&(0,\dots\hat{k_2}\dots,n)=\sigma_{\partial h(
\sigma_\delta)}(0,\dots\hat{k_2}\dots,n)\cr
&\hskip 0.5cm-\varepsilon_{|Deb|}(k_2,0,\dots\hat{k_2}\dots,n)\left(
\sum_{0\leq i<k_0}r_i+\sum_{k_0\leq i<k_1}(r_i-1)\right)\sigma_\delta(0,\dots
\hat{k_2}\dots,n).
\end{align*}
\end{prop}

\vskip 0.3cm
\noindent
{\bf Preuve}\\

On reprend les notations de la proposition pr\'ec\'edente. Le dernier $i$ 
d'ordre $1^+$ dans $\partial\delta$ est $k_2$. On \'ecrit donc:
\begin{align*}
\sigma_{\partial\delta}&=-\sum_{\Delta\in\sigma_\delta}a_\Delta\Big\{
\sum_{\begin{smallmatrix}i<j\cr0\leq i<k_0\cr j=k_2
\end{smallmatrix}}\varepsilon_{|Deb|}(k_2,0,\dots\hat{k_2}\dots,n)
\sum_{\Delta'_{k_2i}\rightarrow^{prop}_i\Delta}\Delta'_{k_2i}
\cr
&\hskip 2cm+\sum_{\begin{smallmatrix}i<j\cr
0\leq i<k_0\cr k_1\leq j<k_2\end{smallmatrix}}\varepsilon_{|Deb|}(j,0,
\dots\hat{\jmath}\dots,n)\sum_{\Delta'_{ji}\rightarrow^{prop}_i\Delta}
\Delta'_{ji}\cr
&\hskip 2cm+\sum_{\begin{smallmatrix}i<j\cr
k_0\leq i<k_1\cr j=k_2\end{smallmatrix}}\varepsilon_{|Deb|}(k_2,0,\dots\hat{k_2}
\dots,n)\Big[\sum_{\Delta'_{k_2i}\rightarrow^{prop}_i\Delta}\Delta'_{k_2i}-
\sum_{\Delta'_{ik_2}\rightarrow^{prop}_i\Delta}\Delta'_{ik_2}\Big]\cr
&\hskip 2cm+\sum_{\begin{smallmatrix}i<j\cr k_0\leq i<k_1\cr k_1\leq j<k_2
\end{smallmatrix}}\varepsilon_{|Deb|}(j,0,\dots\hat{\jmath}\dots,n)\Big[
\sum_{\Delta'_{ji}\rightarrow^{prop}_i\Delta}\Delta'_{ji}-\sum_{\Delta'_{ij}
\rightarrow^{prop}_i\Delta}\Delta'_{ij}\Big]\cr
&\hskip 2cm+\sum_{\begin{smallmatrix}i<j\cr k_1\leq i<k_2\cr j=k_2
\end{smallmatrix}}\varepsilon_{|Deb|}(k_2,0,\dots\hat{k_2}\dots,n)\Big[
\sum_{\Delta'_{k_2i}\rightarrow^{prop}_i\Delta}\Delta'_{k_2i}-
\sum_{\Delta'_{ik_2}\rightarrow^{prop}_i\Delta}\Delta'_{ik_2}\Big]\cr
&\hskip 2cm+\sum_{\begin{smallmatrix}i<j\cr k_1\leq i<k_2\cr k_1\leq j<k_2
\end{smallmatrix}}\varepsilon_{|Deb|}(j,0,\dots\hat{\jmath}\dots,n)\Big[
\sum_{\Delta'_{ji}\rightarrow^{prop}_i\Delta}\Delta'_{ji}-\sum_{\Delta'_{ij}
\rightarrow^{prop}_i\Delta}\Delta'_{ij}\Big]\Big\}.
\end{align*}
Appliquons $h$, on obtient:
\begin{align*}
&h\left(\sigma_{\partial\delta}\right)=-\sum_{\Delta\in\sigma_\delta}a_\Delta
\Big\{\varepsilon_{|Deb|}(k_2,0,\dots\hat{k_2}\dots,n)\cr&\left(
\sum_{0\leq i<k_0}r_i+\sum_{k_0\leq i<k_1}(r_i-1)+\sum_{k_1\leq i<k_2}
(1-1)\right)\Delta(0,\dots\hat{k_2}\dots,n)+~~reste\Big\}.
\end{align*}
Le reste est compos\'e de termes de la forme $h(\Delta'_{ij})$ ou 
$h(\Delta'_{ji})$ avec $i<j<k_2$. Par construction, $h$ \'etant la compression
de la fl\`eche issue du sommet $k_2$ dans $\Delta'_{ji}$ (resp. $\Delta'_{ij}$), 
il r\'ealise une bijection entre les ensembles
$$
\left\{\Delta'_{ji},~~\Delta'_{ji}\rightarrow^{prop}_i\Delta\right\}~~
\stackrel{h}{\longrightarrow}~~\left\{B'_{ji},~~B'_{ji}\rightarrow^{prop}_i
h(\Delta)\right\}
$$
$\Big($resp. entre les ensembles
$$
\left\{\Delta'_{ij},~~\Delta'_{ij}\rightarrow^{prop}_i\Delta\right\}~~
\stackrel{h}{\longrightarrow}~~\left\{B'_{ij},~~B'_{ij}\rightarrow^{prop}_i
h(\Delta)\right\}\Big)
$$
o\`u $B'_{ji}$ (resp. $B'_{ij}$) est un graphe de sommets num\'erot\'es
$(0,\dots\hat{k_2}\dots,n)$ et $h(\Delta)$ a pour sommets $(0,\dots\hat{\jmath}
\dots\hat{k_2}\dots,n)$.

Dans $reste$, chacun des termes correspondant est affect\'e du signe 
$$
\varepsilon_{|Deb|}(j,0,\dots\hat{\jmath}\dots,n).
$$ 
Ce signe co\"\i ncide avec 
le m\^eme signe calcul\'e en supprimant l'indice $k_2$:
$$
\varepsilon_{|Deb(\Delta')|}(j,0,\dots\hat{\jmath}\dots,n)=
\varepsilon_{|Deb(h(\Delta'))|}(j,0,\dots\hat{\jmath}\dots,n).
$$
Donc
$$
reste=-\sigma_{\partial h(\sigma_\delta)}(0,\dots\hat{k_2}\dots,n).
$$


\

\section{Cohomologie des graphes vectoriels}

\

Disons qu'un sommet $i$ d'un graphe vectoriel $\Delta$ est \underbar{simple} si 
$|Fin(i)|\leq1$.

\vskip 0.3cm
\begin{prop}\label{simple} {\rm (Le symbole d'un cocycle ne contient que des 
sommets simples)}

\

Soit $\delta$ une combinaison lin\'eaire sym\'etrique de graphes a\'eriens 
vectoriels. On suppose que $\delta$ est un cocycle ($\partial\delta=0$). Alors
il existe un cobord $\partial\beta$ tel que le symbole 
$\sigma_{\delta-\partial\beta}$ de $\delta-\partial\beta$ ne contient que des 
graphes $\Delta$ dont tous les sommets sont simples.
\end{prop}

\vskip 0.3cm
\noindent
{\bf Preuve}

Puisque $\delta$ est un cocycle, $h(\sigma_{\partial\delta})=0$. Dire que 
$\sigma_\delta$ contient au moins un graphe poss\'edant un sommet non simple, 
c'est dire que dans ${\mathcal O}(\delta)$, $k_1>0$. On a alors:
$$
0=-\varepsilon_{|Deb|}(k_2,0,\dots\hat{k_2}\dots,n)\left(\sum_{0\leq i<k_0}r_i+
\sum_{k_0\leq i<k_1}(r_i-1)\right)\sigma_\delta+
\sigma_{\partial h(\sigma_\delta)}.
$$
Le coefficient $a$ de $\sigma_\delta$ est non nul. Puisque $\sigma_\delta$ n'est 
pas nul, on en d\'eduit que $h(\sigma_\delta)$ n'est pas nul, donc $k_2>k_1$, il 
existe des sommets d'ordre $1^+$ dans les graphes $\Delta$ de $\sigma_\delta$.

Posons $\beta_1=-\frac{1}{a}S(h(\sigma_\delta))$. D'abord par construction, les
graphes apparaissant dans $S(h(\sigma_\delta))$ mais pas dans $h(\sigma_\delta)$
sont d'ordre strictement plus petit que les graphes apparaissant dans 
$h(\sigma_\delta)$, ou:
$$
\sigma_{h(\sigma_\delta)}=\sigma_{S(h(\sigma_\delta))}.
$$
Ensuite on a vu que pour calculer le symbole de $\partial\delta$, on ne 
consid\'erait que les graphes du symbole de $\delta$ donc:
$$
\sigma_{\partial h(\sigma_\delta)}=\sigma_{\partial S(h(\sigma_\delta))}=-
a\sigma_{\partial\beta_1}
$$

Alors,
$$
\sigma_{\partial\beta_1}=\sigma_\delta.
$$ 
Autrement dit $\delta_1=\delta-\partial\beta_1$ a un ordre strictement plus
petit que ${\mathcal O}(\delta)$. Si le symbole de $\delta_1$ a des graphes avec
des sommets non simples, on peut recommencer cette op\'eration. Au bout d'un 
nombre fini d'\'etapes, on arrive sur une combinaison de graphes $\delta-
\partial\beta$ dont le symbole ne contient que des sommets simples:
$$
{\mathcal O}(\delta-\partial\beta)=(1^+,\dots,1^+,1,\dots,1,0,\dots,0,0^-,\dots,
0^-).
$$
Donc tous les graphes de $\delta-\partial\beta$ n'ont eux aussi que des sommets 
simples.

\vskip 0.3cm
\begin{defn} {\rm (Les roues)}

\

Une roue sym\'etrique $R_k$ de longueur $k$ est le sym\'etris\'e de la roue 
simple qui est le graphe $\Delta_k$ ayant $k$ sommets $\{1,\dots,k\}$ et les $k$ 
fl\`eches $\{\overrightarrow{12},\overrightarrow{23},\dots,
\overrightarrow{(k-1)k},\overrightarrow{k1}\}$.
\end{defn}

\vskip  0.3cm
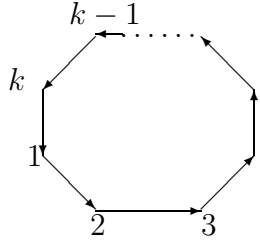
\begin{center}
\begin{picture}(200,100)(100,50)
\put(200,146){\vector(-1,0){10}}
\put(190,145){\vector(-1,-1){20}} \put(170,125){\vector(0,-1){25}}
\put(190,79){\vector(1,0){40}} \put(170,100){\vector(1,-1){20}}
\put(250,100){\vector(0,1){25}}
\put(250,125){\vector(-1,1){20}}\put(199,145){\dots\dots}
\put(230,80){\vector(1,1){20}}\put(164,96){$1$}
\put(188,70){$2$}\put(157,125){$k$}\put(180,150){$k-1$}\put(230,70){$3$}
\end{picture}
\end{center}

\vskip 0.3cm
\begin{lem} {\rm (La cohomologie des roues)}

\

Les roues sym\'etriques de longueur paire sont nulles, $R_{2k}=0$.

Les roues impaires sont des cocycles $\partial R_{2k+1}=0$ qui ne sont pas des 
cobords.
\end{lem}

\vskip 0.3cm
\noindent
{\bf Preuve}

Soit $\Delta_{2k}$ une roue simple de longueur paire. Soit $\sigma$ la 
permutation circulaire $\sigma=(1,2,\dots,2k)$. Alors
$$
\varepsilon(\sigma)=\varepsilon_{|Deb|}(\sigma)=-1
$$
et $\sigma(\Delta_{2k})=\Delta_{2k}$. Donc $R_{2k}=S(\Delta_{2k})=0$.

Par contre les roues $R_{2k+1}$ de longueur impaire ne sont pas nulles. En effet 
si on suppose la dimension de l'espace ${\mathbb R}^d$ assez grande, 
l'op\'erateur
$$
C_{R_{2k+1}}(\alpha_1,\dots,\alpha_{2k+1})=\sum_{\sigma\in S_{2k+1}}
\varepsilon(\sigma)\sum_{1\leq i_1\dots i_{2k+1}\leq d}\partial_{i_{2k+1}}
\alpha^{i_1}_{\sigma(1)}\partial_{i_1}\alpha^{i_2}_{\sigma(2)}\dots
\partial_{i_{2k}}\alpha^{i_{2k+1}}_{\sigma(2k+1)}.
$$
Cet op\'erateur est le cocycle non trivial $\zeta^{(2k+1)}$ de \cite{[DWL]}.
Comme il n'est pas nul, le graphe correspondant n'est pas nul non plus.

Soit $\Delta$ un graphe apparaissant dans la roue sym\'etrique $R_{2k+1}$. Les
graphes $\Delta'_{ij}$ et $\Delta'_{ji}$ apparaissant dans $\partial\Delta$ sont
tous des roues de longueurs $2k+2$ (\`a l'ordre de leur sommets pr\`es). Comme
$\partial R_{2k+1}$ est sym\'etrique, on a donc $\partial R_{2k+1}=0$.

Supposons que $R_{2k+1}=\partial\beta$. Il est clair que $\beta$ ne contient que
des graphes ayant $2k$ sommets tous d'ordres $1^+$. Donc $\beta$ est une 
combinaison lin\'eaire de roues de longueur paire, $\beta=0$, ce qui est 
impossible. $R_{2k+1}$ n'est pas un cobord.

\vskip 0.3cm
\begin{lem} {\rm (Les graphes \`a roues)}

\

Si $\Delta$ est un graphe dont toutes les composantes connexes sont des roues
de longueurs impaires, alors le sym\'etris\'e de $\Delta$ est nul si deux
composantes connexes ont la m\^eme longueur, sinon $S(\Delta)$ est au signe
pr\`es le sym\'etris\'e du graphe:
$$
\Delta_{k_1}\wedge\Delta_{k_2}\wedge\dots\wedge\Delta_{k_p}
$$
dont les fl\`eches sont:
\begin{align*}
\Big\{&\overrightarrow{12},\dots,\overrightarrow{(k_1-1)k_1},
\overrightarrow{k_11},\cr
&\overrightarrow{(k_1+1)(k_1+2)},\dots,\overrightarrow{(k_1+k_2-1)(k_1+k_2)},
\overrightarrow{(k_1+k_2)(k_1+1)},\dots,\cr
&\overrightarrow{(k_1+\dots+k_{p-1}+1)(k_1+\dots+k_{p-1}+2)},\dots,
\overrightarrow{(k_1+\dots+k_p)(k_1+\dots+k_{p-1}+1)}\Big\}
\end{align*}
et $k_1<k_2<\dots<k_p$.
\end{lem}

\vskip 0.3cm
\noindent
{\bf Preuve}

Si deux composantes connexes de $\Delta$ sont des roues de m\^eme longueur 
$2k+1$, de sommets num\'erot\'es $\{i_1,\dots,i_{2k+1},j_1,\dots,j_{2k+1}\}$,
alors la permutation
$$
\sigma=(i_1,j_1)\dots(i_{2k+1},j_{2k+1})
$$ 
est impaire et laisse $\Delta$ invariant, donc le sym\'etris\'e de $\Delta$ est 
nul.

Il est clair que si les composantes connexes de $\Delta$ sont toutes des roues
de longueurs impaires diff\'erentes, il existe une permutation des sommets de
$\Delta$ qui le transforme en $\Delta_{k_1}\wedge\Delta_{k_2}\wedge\dots\wedge
\Delta_{k_p}$ avec $k_1<k_2<\dots<k_p$. On notera abusivement:
$$
S(\Delta)=R_{k_1}\wedge R_{k_2}\wedge\dots\wedge R_{k_p}.
$$

\vskip 0.3cm
\begin{lem} {\rm (La cohomologie des graphes \`a roues)}

\

Tous les graphes \`a roues $R=\bigwedge_{i=1}^pR_{2k_i+1}$ sont des cocycles qui
sont lin\'eairement ind\'edependants dans l'espace de cohomologie.
\end{lem}

\vskip 0.3cm
\noindent
{\bf Preuve}

D'abord chacun de ces graphes est un cocycle, ensuite si une combinaison 
lin\'eaire de ces graphes sym\'etriques est un cobord ($\sum_ja_jR_j=\partial
\beta$), en se pla\c cant sur un espace ${\mathbb R}^d$ de dimension assez 
grande, on obtient une combinaison lin\'eaire d'op\'erateurs
$$
\sum_ja_j\left(\bigwedge_{i=1}^{p_j}\zeta^{(2k_i+1)}\right)
$$
qui est un cobord pour la cohomologie de Chevalley de l'alg\`ebre de Lie des 
champs de vecteurs sur ${\mathbb R}^d$ \`a valeur dans les fonctions (les 
0-formes). Ceci n'est possible que si chaque $a_j$ est nul d'apr\`es 
\cite{[DWL]}.

\vskip 0.3cm
\begin{defn} {\rm (Les lignes)}

\

La ligne sym\'etrique $L_\ell$ de longueur $\ell$ est le sym\'etris\'e de la
ligne simple $\Delta_\ell$, graphe de sommets $\{1,\dots,\ell+1\}$ et de 
fl\`eches $\{\overrightarrow{21},\overrightarrow{32},\dots,
\overrightarrow{(\ell+1)\ell}\}$. Par convention, la ligne $L_0$ est le graphe 
\`a un sommet $\{1\}$ et sans fl\`eche.
\end{defn}

\vskip 0.3cm
\begin{lem} {\rm (Cohomologie des lignes)}

\

Les lignes sym\'etriques de longueur impaire $L_{2\ell+1}$ sont des cobords.

Le sym\'etris\'e d'un graphe contenant deux lignes de m\^eme longueur impaire
est nul.

Les lignes sym\'etriques de longueur paire $L_{2\ell}$ ne sont pas des cocycles, 
plus pr\'eci\-s\'e\-ment, on a:
$$
\partial L_{2\ell}=L_{2\ell+1}.
$$
\end{lem}

\vskip 0.3cm
\noindent
{\bf Preuve}

Prenons une ligne simple $\Delta_\ell$. Si $\ell>0$, dans le calcul de $\partial
\Delta_\ell$, on \'eclate les $\ell-1$ sommets qui ne sont pas les 
extr\'emit\'es de la ligne:
$$
\partial\Delta_\ell=-\sum_{i\neq j}\varepsilon_{|Deb|}(j,0,\dots\hat{\jmath}
\dots,\ell+1)\sum_{\Delta'_{ji}\rightarrow^{prop}_i\Delta_\ell}\Delta'_{ji}.
$$
Pour chaque couple $i\neq j$, il y a un seul graphe $\Delta'_{ji}$. A l'ordre 
des sommets pr\`es, on obtient \`a chaque fois une ligne de longueur $\ell+1$.

On cherche donc dans le cobord de $L_\ell$ les lignes simples $\Delta_{\ell+1}$.
Ces lignes ne peuvent provenir que de l'\'eclatement du sommet $i$ de la ligne
simple $\Delta_\ell$ avec les sommets $\{0,\dots\hat{\imath+1}\dots,\ell+1\}$ et 
pour $j=i+1$. On obtient donc
\begin{align*}
-\left(\sum_{i=1}^\ell\varepsilon_{|Deb|}(i+1,0,\dots\widehat{\imath+1}\dots,
\ell+1)\right)\Delta_{\ell+1}&=-\left(\sum_{i=1}^\ell(-1)^i\right)
\Delta_{\ell+1}
\cr
&=\left\{\begin{array}{cccc}0&~~\hbox{si}&~~\ell&~~\hbox{est pair},\cr 
\Delta_{\ell+1}&~~\hbox{si}&~~\ell&~~\hbox{est impair}.\end{array}\right.
\end{align*}
En sym\'etrisant, on obtient:
$$
\partial L_{2\ell+1}=0,\quad \partial L_{2\ell}=L_{2\ell+1}.
$$

Si $\Delta$ contient deux composantes connexes qui sont des lignes de m\^eme
longueur impaire $2\ell+1$, de sommets num\'erot\'es $i_1,\dots,i_{2\ell+2}$ et
$j_1,\dots,j_{2\ell+2}$, la permutation
$$
\sigma=(i_1,j_1)\dots(i_{2\ell+2},j_{2\ell+2})
$$
est telle que:
$$
\varepsilon_{|Deb|}(\sigma)=-1\quad\hbox{et}\quad \sigma(\Delta)=\Delta.
$$
Donc le sym\'etris\'e $S(\Delta)$ de $\Delta$ est nul.

\vskip 0.3cm
\begin{lem} {\sl (Les graphes \`a lignes)}

\

Si $\Delta$ est un graphe dont toutes les composantes connexes sont des lignes,
si le sym\'etris\'e $S(\Delta)$ de $\Delta$ n'est pas nul, c'est au signe pr\`es 
le sym\'etris\'e du graphe:
$$
\left(\Delta_0^{k_0}\Delta_2^{k_1}\dots\Delta_{2\ell}^{k_\ell}\right)
\Delta_{2\ell_1+1}\wedge\Delta_{2\ell_2+1}\wedge\dots\wedge\Delta_{2\ell_q+1}
\qquad(\ell_1<\ell_2<\dots<\ell_q).
$$
Les sommets de ce graphe sont rang\'es dan l'ordre naturel (par exemple
$\Delta_{2\ell}^{k_\ell}$ d\'esigne un graphe ayant $k_\ell(2\ell+1)$ sommets
form\'e d'une union de $k_\ell$ lignes de longueur $2\ell$, en num\'erotant 
d'abord les sommets de la premi\`ere ligne, puis ceux de la seconde, etc\dots)
\end{lem}

\vskip 0.3cm
\noindent
{\bf Preuve}

Si le sym\'etris\'e de $\Delta$ n'est pas nul, pour chaque longueur impaire, il 
y a au plus une composante connexe de cette longueur. Les composantes de 
longueur paire (y compris 0) peuvent appara\^\i tre plusieurs fois. Modulo une
permutation des sommets, $\Delta$ est donc bien un produit de lignes comme
annonc\'e.

Notons abusivement le sym\'etris\'e de ce graphe 
$$
\left(L_0^{k_0}L_2^{k_1}\dots L_{2\ell}^{k_\ell}\right)
L_{2\ell_1+1}\wedge L_{2\ell_2+1}\wedge\dots\wedge L_{2\ell_q+1}
$$

\vskip 0.3cm
\begin{lem} {\rm (La cohomologie des graphes \`a lignes)}

\

Consid\'erons le graphe \`a lignes:
$$
\left(\prod_{i=0}^\ell L_{2i}^{k_i}\right)\left(\bigwedge_{j=1}^q
L_{2\ell_j+1}\right).
$$
Si pour chaque $i$, $k_i\neq0$ implique qu'il existe $j$ tel que 
$i=\ell_j$, alors
$$
\partial\left(\prod_{i=0}^\ell L_{2i}^{k_i}\right)\left(\bigwedge_{j=1}^q
L_{2\ell_j+1}\right)=0.
$$
Sinon
$$
\partial\left(\prod_{i=0}^\ell L_{2i}^{k_i}\right)\left(\bigwedge_{j=1}^q
L_{2\ell_j+1}\right)=\sum_{r=0}^\ell k_r \left(\prod_{i\neq r}L_{2i}^{k_i}
\right)L_{2r}^{k_r-1}L_{2r+1}\wedge\left(\bigwedge_{j=1}^q
L_{2\ell_j+1}\right)\neq0.
$$

La cohomologie des graphes \`a lignes est triviale.
\end{lem}

\vskip 0.3cm
\noindent
{\bf Preuve}

Il r\'esulte des calculs pr\'ec\'edents et du fait que les lignes de longueur 
paire peuvent permuter avec toutes les autres lignes sans changement de signe 
que:
$$
\partial\left(\prod_{i=0}^\ell L_{2i}^{k_i}\right)\left(\bigwedge_{j=1}^q
L_{2\ell_j+1}\right)=\sum_{r=0}^\ell k_r \left(\prod_{i\neq r}L_{2i}^{k_i}
\right)L_{2r}^{k_r-1}L_{2r+1}\wedge\left(\bigwedge_{j=1}^q
L_{2\ell_j+1}\right).
$$
Donc si pour chaque $i$, $k_i\neq0$ implique qu'il existe $j$ tel que 
$i=\ell_j$, alors
$$
\partial\left(\prod_{i=0}^\ell L_{2i}^{k_i}\right)\left(\bigwedge_{j=1}^q
L_{2\ell_j+1}\right)=0.
$$

S'il existe $i$ tel que $k_i>0$ et il n'y a pas de $j$ tel que $i=\ell_j$, le
graphe n'est pas un cocycle car si $r\neq s$, 
$$
\left(\prod_{i\neq r}L_{2i}^{k_i}\right)L_{2r}^{k_r-1}L_{2r+1}\wedge\left(
\bigwedge_{j=1}^qL_{2\ell_j+1}\right)\neq\left(\prod_{i\neq s}L_{2i}^{k_i}
\right)L_{2s}^{k_s-1}L_{2s+1}\wedge\left(\bigwedge_{j=1}^qL_{2\ell_j+1}\right).
$$
En particulier le graphe \`a lignes simples
$$
\left(\prod_{i=0}^{\ell-1}\Delta_{2i}^{k_i}\right)\Delta_{2\ell}^{k_\ell-1}
\Delta_{2\ell+1}\wedge\left(\bigwedge_{j=1}^q\Delta_{2\ell_j+1}\right)
$$
n'appara\^\i t qu'une seule fois. Le second membre n'est pas nul.

Si un graphe \`a lignes
$$
L=\left(\prod_{i=1}^p L_{2i}^{k_i}\right)\left(\bigwedge_{j=1}^qL_{2\ell_j+1}
\right)$$
est un cocycle, $2\ell_q+1$ est la plus grande longueur des graphes 
apparaissant dans $L$ et on a si $p=\ell_q$
$$
L=(-1)^{q-1}\partial\frac{1}{k_p+1}\left(\prod_{i=1}^{p-1} L_{2i}^{k_i}\right)
L_{2p}^{k_p+1}\left(\bigwedge_{j=1}^{q-1}L_{2\ell_j+1}\right)
$$
et si $p<\ell_q$
$$
L=(-1)^{q-1}\partial\left(\prod_{i=1}^p L_{2i}^{k_i}\right)L_{2\ell_q}\left(
\bigwedge_{j=1}^{q-1}L_{2\ell_j+1}\right).
$$
La cohomologie des graphes \`a lignes est donc toujours triviale.

\vskip 0.3cm
\begin{thm} {\rm (La cohomologie de Chevalley des graphes)}

\

La cohomologie de Chevalley des graphes vectoriels est donn\'ee par les roues de
longueur impaire. Plus pr\'ecis\'ement, pour tout $n$, une base de $H^n$ est 
donn\'ee par
$$
\left\{R_{2k_1+1}\wedge R_{2k_2+1}\wedge\dots\wedge R_{2k_p+1},\quad\hbox{avec}
~~k_1<k_2<\dots<k_p,~~\sum_{i=1}^p(2k_i+1)=n\right\}.
$$
\end{thm}

\vskip 0.3cm
\noindent
{\bf Preuve}

D'apr\`es la proposition \ref{simple}, tout cocycle $\delta$ est cohomologue \`a 
un cocycle $\delta-\partial\beta$ dont le symbole ne contient que des graphes 
avec des sommets simples. On suppose maintenant que $\delta$ a cette 
propri\'et\'e. Chaque graphe de ce symbole est donc une union de composantes 
connexes qui sont soit des roues simples de longueur impaire soit des lignes 
simples. Le nombre de lignes est d'ailleurs fix\'e, \'egal au nombre de $0^-$.

On d\'efinit un nouvel ordre sur les graphes de ce type en posant:
$$
{\mathcal O}'(\Delta)=(\ell_1,\dots,\ell_p,r_1,\dots,r_q)
$$
o\`u les $\ell_i$ sont les longueurs des lignes et les $r_j$ celles des roues.
On range ces symboles en posant $\ell_i>r_j$ pour tout $i$ et $j$, on range ces 
ordres de graphes par l'ordre lexicographique sur les symboles.

L'ordre de $\delta$ est
$$
{\mathcal O}'(\delta)=Max\{{\mathcal O}'(\Delta),~~\Delta~~\hbox{apparaissant 
dans}~~\sigma_\delta\}.
$$
On a donc ${\mathcal O}'(\delta)=(\ell_1,\dots,\ell_p,r_1,\dots,r_q)$ avec 
$\ell_1\geq\ell_2,\dots,r_1>r_2>\dots$, les seuls $\ell_i$ r\'ep\'et\'es sont 
pairs, les $r_j$ sont tous impairs. Le nouveau symbole de $\delta$ est 
$$
\sigma'_\delta=\sum_{\begin{smallmatrix}\Delta\cr
{\mathcal O}'(\Delta)={\mathcal O}'(\delta)\end{smallmatrix}}a_\Delta\Delta.
$$\\
Supposons qu'il existe un indice $i$ tel que $\ell_i$ est pair et $\ell_i+1$ 
n'appartient pas \`a $\{\ell_1,\dots,\ell_{i-1}\}$. Appelons $i_0$ le premier 
indice pour lequel ceci se produit.

Si $\Delta$ appara\^\i t dans $\delta$ mais pas dans $S(\sigma_\delta)$, on a vu
que:
$${\mathcal O}(\partial\Delta)<{\mathcal O}(\delta)\oplus1^+.
$$
D'apr\`es les lemmes pr\'ec\'edents, si $\Delta$ appara\^\i t dans 
$S(\sigma_\delta)$ mais pas dans $S(\sigma'_\delta)$, alors le nouvel ordre 
${\mathcal O}'(S(\partial\Delta))$ est strictement plus petit que
$$
{\mathcal O}'(\delta)\oplus [i_0]=(\ell_1,\dots,
\underbrace{\ell_{i_0}+1}_{(i_0)},\dots,\ell_p,r_1,\dots,r_q).
$$
En effet, pour calculer $\partial\Delta$, on proc\`ede \`a des \'eclatements
soit dans les lignes soit dans les roues. On ne retient que les \'eclatements 
qui ne disparaissent pas par sym\'etrisation.\\

Enfin si $\Delta$ appara\^\i t dans $S(\sigma'_\delta)$, le nouvel ordre de 
$S(\partial\Delta)$ est toujours inf\'erieur ou \'egal \`a ${\mathcal O}'
(\delta)\oplus[i_0]$, car on utilise les lemmes pr\'ec\'edents, en proc\'edant
\`a des \'eclatements successivement des sommets des lignes et des sommets des 
roues.

On a donc avec nos notations habituelles:
$$
S(\sigma'_\delta)=a\left(\prod_{\ell_i~pair}L_{\ell_i}\right)\left(
\bigwedge_{\ell_i~impair}L_{\ell_i}\right)\wedge\left(\bigwedge_{j=1}^pR_j\right)
$$
et
$$
S(\sigma'_{\partial\delta})=a\left(\prod_{\begin{smallmatrix}\ell_i~pair\cr
i\neq i_0\end{smallmatrix}}L_{\ell_i}\right)L_{\ell_{i_0}+1}\wedge\left(
\bigwedge_{\ell_i~impair}L_{\ell_i}\right)\wedge\left(\bigwedge_{j=1}^pR_j
\right)\neq0.
$$
Ceci est donc impossible, il n'existe pas de $i$ tel que $\ell_i$ soit pair et
$\ell_i+1$ n'appartient pas \`a $\{\ell_1,\dots,\ell_{i-1}\}$.\\

On en d\'eduit en particulier que s'il y a des lignes, $\ell_1$ est impair et en 
retranchant \`a $\delta$ le cobord de
$$
\beta=a\left(\prod_{\ell_i~pair}L_{\ell_i}\right)L_{\ell_1-1}\left(
\bigwedge_{\begin{smallmatrix}\ell_i~impair\cr i>1\end{smallmatrix}}L_{\ell_i}
\right)\wedge\left(\bigwedge_{j=1}^pR_j\right)
$$
on obtient un graphe $\delta-\partial\beta$ dont le nouvel ordre est strictement
plus petit que celui de $\delta$. En r\'ep\'etant cette op\'eration, on se
ram\`ene au cas o\`u le symbole de $\delta$ ne contient que des graphes \`a 
roues impaires. Alors $S(\sigma_\delta)$ est un cocycle et l'ordre de $\delta-
S(\sigma_\delta)$ est strictement plus petit.\\

On recommence cette op\'eration jusqu'\`a annulation de l'ordre des cocycles
cons\-truits. On obtient :
$$
\delta=\partial\beta+\sum_{r_1,\dots,r_q}a_{r_1\dots r_q}R_{r_1}\wedge\dots
\wedge R_{r_q}.
$$
La cohomologie est donc engendr\'ee par les graphes \`a roues de longueurs 
impaires. On a vu que ces graphes sont lin\'eairement ind\'ependants.

\vfill
\eject

\

\end{document}